\documentclass [11pt]{amsart}

\usepackage{graphicx}
\usepackage{amsmath}
\usepackage{amssymb}

\newcommand{\nin}{\notin}
\newcommand{\bdry}{\partial}

\renewcommand{\Tilde}{\widetilde}

\newcommand{\cardioid}{\heartsuit}

\newcommand{\pf}[1]{\noindent{\bf Proof#1:\  }}
\newcommand{\defn}[1]{\medskip\noindent{\bf Definition#1:\  }}
\newcommand{\obs}[1]{\noindent{\bf Observation#1:\ }}
\newcommand{\note}[1]{\noindent{\bf Note#1:\  }}

\newcommand{\QED}{\hfill$\square$\medskip}

\newtheorem{thm}{Theorem}[section]
\newtheorem{prop}[thm]{Proposition}
\newtheorem{lemma}[thm]{Lemma}
\newtheorem{corol}[thm]{Corollary}

\begin{document}

\title{Quadratic polynomials and combinatorics of the principal nest}
\author{Rodrigo A. P\'erez}
\address{Department of Mathematics, Cornell University, Ithaca, NY
14853. USA.}
\email{rperez@math.cornell.edu}
\thanks{Research supported by an NSF Postdoctoral Fellowship in the
Mathematical Sciences, grant DMS-0202519.}

\begin{abstract}
  {\it The definition of principal nest is supplemented with a system of
  frames that make possible the classification of combinatorial types for
  every level of the nest. As a consequence, we give necessary and sufficient
  conditions for the admissibility of a type and prove that given a sequence
  of non-renormalizable finite admissible types, there is a quadratic
  polynomial whose nest realizes the sequence.}
\end{abstract}

\maketitle

\renewcommand{\theequation}{\arabic{section}.\arabic{equation}}

\section{Introduction}\label{sect:Intro}
\setcounter{equation}{0}

We will study the combinatorial behavior of the dynamics for quadratic
polynomials with (non-periodic) recurrent critical orbit; these are the maps
that have a well defined {\it principal nest}. 

In \cite{L_nest}, M. Lyubich developed the principal nest as a tool to provide
some examples of infinitely renormalizable parameters at which the Mandelbrot
set is locally connected. The nest consists of a subsequence of central puzzle
pieces, each determined by the first return of the critical orbit to the
preceding nest piece.

As described in Section 2, the principal nest may include non-central pieces
at some levels. Each piece $V$ of the nest has a first return map onto the
central piece of previous level that contains $V$.

When the polynomial is real, the lateral pieces of the nest can only be
located to the left or right of the central piece. This information, together
with the sign of the derivative of the first return maps, is enough to provide
a complete classification of real nest types (see
\cite{L_attractor}). However, in the complex case, lateral pieces may ``hang''
from different branches of the Julia set. We exploit this underlying structure
to construct a {\it frame system} that encodes the configuration of the
nest. This allows us to describe the possible itineraries of the critical
orbit as it visits different levels.

Our main classification result is the following:

\medskip\noindent{\bf Theorem: }{\it Any infinite sequence of finite, weak
  combinatorial types is realized in the quadratic family, as long as the
  types satisfy the admissibility condition at every level. The set of
  parameters that display this sequence of types can be described as the
  residual intersection in an infinite family of sequences of nested
  parapieces.
\medskip}

We illustrate the applicability of frames with a description of maximal
hyperbolic components of the Mandelbrot set, and with the construction of
complex analogues of the {\it rotation-like maps} of
\cite{B_odometers}. Further applications, including a classification of
complex quadratic Fibonacci maps, are contained in \cite{2nd_part}.

\subsection{Background and organization}
The concept of a puzzle partition was introduced in \cite{BH_cubics1} and
\cite{BH_cubics2} to study the topology of cubic Julia sets as a function of
the critical points. In the late 80's, J.-C. Yoccoz implemented the puzzle in
the setting of quadratic polynomials, in order to prove the MLC conjecture for
the case of finitely renormalizable parameters (see \cite{Tableaux}). The idea
of the puzzle construction is to show that the pieces around the critical
point become arbitrarily small, thus providing a system of neighborhoods that
satisfy the local connectivity condition. For Yoccoz's puzzle, this is done by
showing that the moduli of annuli between consecutive pieces generate a
divergent series. In the case of the principal nest, the moduli between
consecutive nest pieces increase in an essentially linear fashion. The
principal nest technique underlies Lyubich's proofs of the
Feigenbaum-Collet-Tresser conjecture and the theorem on the measure-theoretic
attractor. \\

In order to fix notation, we introduce basic notions of Complex Dynamics in
Section \ref{sect:Basics}. In particular, we describe the puzzle construction
of Yoccoz and the principal nest following Lyubich.

In Section \ref{sect:Frames} we define the frame associated to a nest. The
construction requires particular care at the initial steps in order to ensure
that nest levels and frame levels go hand by hand. Then we specify a labeling
of frame cells and produce a language to describe admissible combinatorial
types. Our main result (Theorem \ref{thm:Frame_Renormalizations} and Corollary
\ref{corol:Infinite_Type}) is stated and proved there.

Section \ref{sect:Examples} illustrates the use of our construction with two
examples; a classification of {\it maximal} hyperbolic components of the
Mandelbrot set according to the combinatorial type of their nests, and an
extension of the family of {\it rotation-like maps} described in
\cite{B_odometers}.

A brief summary of holomorphic motions is included in an appendix.

\subsection{Acknowledgments}
This work contains results from my dissertation. Many thanks are due to my
advisors John Milnor and Mikhail Lyubich for their generous support during the
preparation of the Thesis. I would also like to thank John Smillie for
suggestions to improve the presentation. Finally, some of the pictures were
created with the PC program {\tt mandel.exe} by Wolf Jung \cite{J}.

\section{Basics in Complex Dynamics}\label{sect:Basics}
\setcounter{equation}{0}

\subsection{Basic notions}
In order to fix notation, let us start by defining the basic notions of
complex dynamics that will be used; we refer the reader to \cite{DH_orsay} and
\cite{M_book} for details on this introductory material.

\begin{center}\begin{figure}[h]
  \includegraphics{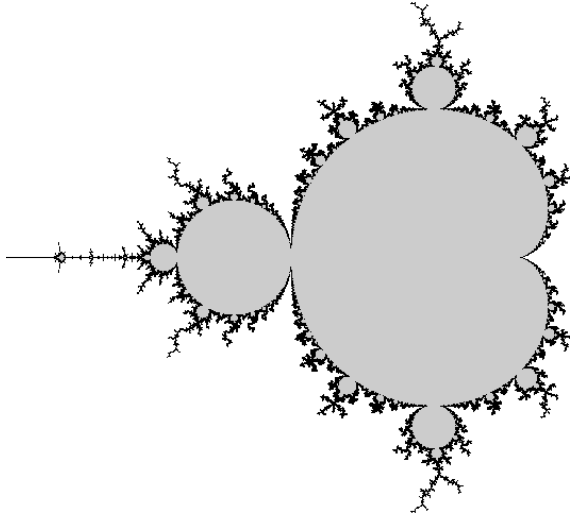}
  \caption[The Mandelbrot set]{\label{fig:Mandel} \it The Mandelbrot set.}
\end{figure}\end{center}

We focus attention on the {\it quadratic family} $\mathcal{Q} := \big\{ f_c:z
\mapsto z^2 + c \mid c \in \mathbb{C} \big\}$. For every $c$, the compact sets
$K_c := \big\{ z \mid \text{the sequence } \{f_c^{\circ n}(z)\} \text{ is
bounded} \big\}$ and $J_c := \bdry K_c$ are called the {\bf filled Julia set}
and {\bf Julia set} respectively. Depending on whether the orbit of the
critical point 0 is bounded or not, $J_c$ and $K_c$ are connected or totally
disconnected. The {\bf Mandelbrot set} is defined as $M := \big\{c \mid c
\in K_c \big\}$; that is, the set of parameters with bounded critical orbit;
see Figure \ref{fig:Mandel}.

A component of $\text{int} \,M$ that contains a superattracting parameter will
be called a {\bf hyperbolic component}\footnote{Though, of course, it is
conjectured that all interior components are hyperbolic.}. The boundary of a
hyperbolic component can either be real analytic, or fail to be so at one cusp
point. The later kind are called {\bf primitive} components. In particular,
the hyperbolic component $\cardioid$ associated to $z \mapsto z^2$ is bounded
by a cardioid known as the {\bf main cardioid}.

$M$ contains infinitely many small homeomorphic copies of itself, accumulating
densely around $\bdry M$. In fact, every hyperbolic component $H$ other than
the main one is the base of one such small copy $M'$. $H$ is called {\bf
prime} if it is not contained in any other small copy. To simplify later
statements, prime components are further subdivided in {\bf immediate}
(non-primitive components that share a boundary point with $\cardioid$) and
{\bf maximal} (primitive components away from $\bdry\cardioid$).

\subsection{External rays, wakes and limbs}
  \label{subsect:External_Rays_n_Wakes} Since $f_c^{-1}(\infty) =\{\infty\}$,
the point $\infty$ is a fixed critical point and a result of B\"ottcher yields
a change of coordinates that conjugates $f_c$ to $z \mapsto z^2$ in a
neighborhood of $\infty$. With the requirement that the derivative at $\infty$
is 1, this conjugating map is denoted $\varphi_c:N_c \longrightarrow
\overline{\mathbb{C}} \setminus \overline{\mathbb{D}_R}$, where $\mathbb{D}_R$
is the disk of radius $R \geq 1$ and $N_c$ is the maximal domain of
unimodality for $\varphi_c$. It can be shown that $N_c = \overline{\mathbb{C}}
\setminus K_c$ and $R=1$ whenever $c \in M$. Otherwise, $N_c$ is the exterior
of a figure 8 curve that is real analytic and symmetric with respect to 0. In
this case, $R>1$ and $K_c$ is contained in the two bounded regions determined
by the 8 curve.

Consider the system of radial lines and concentric circles in $\mathbb{C}
\setminus \mathbb{D}_R$ that characterizes polar coordinates. The pull back of
these curves by $\varphi_c$, creates a collection of {\bf external rays}
$r_{\theta}$ $\big( \theta \in [0, 1) \big)$ and {\bf equipotential curves}
$e_s$ \big(here $s \in (R,\infty)$ is called the {\bf radius} of $e_s$\big) on
$N_c$. These form two orthogonal foliations that behave nicely under dynamics:
$f_c(r_{\theta}) = r_{2\theta}$, $f_c(e_s) = e_{(s^2)}$. When $c \in M$, we
say that a ray $r_{\theta}$ {\bf lands} at $z \in J_c$ if $z$ is the only
point of accumulation of $r_{\theta}$ on $J_c$. \\

A similar coordinate system exists around the Mandelbrot set. For $c \nin M$,
we define the map
\begin{equation}\label{eqn:Phi_M}
  \Phi_M(c) := \varphi_c(c).
\end{equation}

In \cite{DH_orsay} it is shown that $\Phi_M:\mathbb{\overline{C}} \setminus M
\longrightarrow \mathbb{\overline{C}} \setminus \overline{\mathbb{D}}$ is a
conformal homeomorphism tangent to the identity at $\infty$. This yields
connectivity of $M$ and allows us to define {\bf parametric external rays} and
{\bf parametric equipotentials} as in the dynamical case. Since there is
little risk of confusion, we will use the same notation ($r_{\theta}, e_s$) to
denote these curves and say that a parametric ray lands at a point $c \in
\bdry M$ if $c$ is the only point of accumulation of the ray on $M$.

For the rest of this work, all rays considered, whether in dynamical or
parameter plane, will have rational angles. These are enough to work out our
combinatorial constructions and satisfy rather neat properties.

\begin{prop} (\cite{M_book}, ch.18) Both in the parametric and the dynamical
  situations, if $\theta \in \mathbb{Q}$ the external ray $r_{\theta}$
  lands. In the dynamical case, the landing point is (pre-)periodic with the
  period and preperiod determined by the binary expansion of $\theta$. A point
  in $J_c$ (respectively $\bdry M$) can be the landing point of at most, a
  finite number of rays (respectively parametric rays). If this number is
  larger than 1, each component of the plane split by the landing rays will
  intersect $J_c$ (respectively $\bdry M$).
\end{prop}

Unless $c = \frac14$, $f_c$ has two distinct fixed points. If $c \in M$,
these can be distinguished since one of them is always the landing point of
the ray $r_0$. We call this fixed point $\beta$. The second fixed point is
called $\alpha$ and can be attracting, indifferent or repelling, depending on
whether the parameter $c$ belongs to $\cardioid$, $\bdry\cardioid$, or
$\mathbb{C} \setminus \overline{\cardioid}$. The map $\psi_0: \cardioid
\longrightarrow \mathbb{D}$ given by $c \mapsto f'_c(\alpha_c)$ is the Riemann
map of $\cardioid$ normalized by $\psi_0(0)=0$ and $\psi'_0(0)>0$. Since the
cardioid is a real analytic curve except at $\frac14$, $\psi_0$ extends to
$\overline{\cardioid}$.

The fixed point $\alpha$ is parabolic exactly at parameters $c_{\eta} \in
\bdry \cardioid$ of the form $c_{\eta} = \psi_0^{-1} \left( e^{2\pi i \eta}
\right)$ where $\eta \in \mathbb{Q} \cap [0,1)$. If $\eta \neq 0$, $c_{\eta}$
is the landing point of two parametric rays $r_{t^-(\eta)}$ and
$r_{t^+(\eta)}$.

\defn{} The closure of the component of $\mathbb{C} \setminus \left(
  r_{t^-(\eta)} \cup c_{\eta} \cup r_{t^+(\eta)} \right)$ that does not
  contain $\cardioid$ is called the {\bf $\eta$-wake} of $M$ and is denoted
  $W_{\eta}$. The {\bf $\eta$-limb} is defined as $L_{\eta} = M \cap
  W_{\eta}$.
\medskip

\defn{} Say that $\eta=\frac pq$, written in lowest terms. Then
  $\mathcal{P}\big( \frac pq \big)$ will denote the unique set of angles whose
  behavior under doubling is a cyclic permutation with combinatorial rotation
  number $\frac pq$.
\medskip

If $\mathcal{P}\big( \frac pq \big) = \{t_1,\ldots , t_q \}$, then for any
parameter $c \in L_{p/q}$ the corresponding point $\alpha$ splits $K_c$ in $q$
parts, separated by the $q$ rays $\{ r_{t_1}, \ldots , r_{t_q} \}$ landing at
$\alpha$. The two rays whose angles span the shortest arc separate the
critical point 0 from the critical value $c$; these two angles turn out to be
$t^-(\frac pq)$ and $t^+(\frac pq)$.

\subsection{Yoccoz puzzles}
The Yoccoz {\bf puzzle} is well defined for parameters $c \in L_{p/q}$ for any
any $\frac pq \in \mathbb{Q} \cap [0,1)$ with $(p,q)=1$. If 0 is not a
preimage of $\alpha$, the puzzle is defined at infinitely many depths and we
will restrict attention to these parameters. Since we describe properties of a
general parameter, we will omit the subscript and write $f$ instead of $f_c$,
$K$ instead of $K_c$ and so on.

Let us fix the neighborhood $U$ of $K$ bounded by the equipotential of radius
2. The rays that land at $\alpha$ determine a partition of $U \setminus \{
r_{t_1}, \ldots , r_{t_q} \}$ in $q$ connected components.  We will call the
closures $Y_0^{(0)}, Y_1^{(0)}, \ldots, Y_{q-1}^{(0)}$ of these components,
{\bf puzzle pieces} of depth $0$. At this stage the labeling is chosen so that
$0 \in Y_0^{(0)}$ and $f \left( K \cap Y_j^{(0)} \right) = K \cap
Y_{j+1}^{(0)}$; where the subindices are understood as residues modulo $q$. In
particular, $Y_1^{(0)}$ contains the critical value $c$ and the angles of its
bounding rays are $t^-(\frac pq),t^+(\frac pq)$.

The puzzle pieces $Y_i^{(n)}$ of higher depths are recursively defined as the
closures of every connected component in $f^{\circ (-n)} \left( \bigcup
\text{int}\, Y_j^{(0)} \right)$; see Figure \ref{fig:Puzzle_Graphs}. At each
depth $n$, there is a unique piece which contains the critical point and we
will always choose the indices so that $0 \in Y_0^{(n)}$.

We will denote by $P_n$ the collection of pieces of level $n$. The resulting
family $\mathcal{Y}_c := \{ P_0, P_1, \ldots \}$ of puzzle pieces of all
depths, has the following two properties:

{\renewcommand{\labelenumi}{{\bf P\arabic{enumi}}}
\renewcommand{\theenumi}{\labelenumi}
\begin{enumerate}
  \item \label{piece_relns} Any two puzzle pieces either are nested (with the
        piece of higher depth contained in the piece of lower depth), or have
        disjoint interiors.

  \item The image of any piece $Y_j^{(n)}$ $(n \geq 1)$ is a piece
        $Y_i^{(n-1)}$ of the previous depth $n-1$. The restricted map
        $f:\text{int}\, Y_j^{(n)} \longrightarrow \text{int}\, Y_i^{(n-1)}$
        is a 2 to 1 branched covering or a conformal homeomorphism, depending
        on whether $j=0$ or not.
\end{enumerate}}

These properties characterize $\mathcal{Y}_c$ as a Markov family, endowing the
puzzle partition with dynamical meaning.

Note that the collection of ray angles at depth $n$ consists of all
$n$-preimages of $\{ r_{t_1}, \ldots , r_{t_q} \}$ under angle doubling. The
union of all pieces of depth $n$ is the region enclosed by the equipotential
$e_{\left( 2^{2^{-n}} \right)}$. Note also that every piece $Y$ of depth $n$
is the $n^{\text{th}}$ preimage of some piece of level 0. By further
iteration, $Y$ will map onto a region determined by the same rays as
$Y_0^{(0)}$ and a possibly larger equipotential. This provides a 1 to 1
correspondence between puzzle pieces and preimages of 0. The distinguished
point inside each piece is called the {\bf center} of the piece.

\subsection{Adjacency Graphs}\label{subsect:Graphs}
Given a set of puzzle pieces $P \subset P_n$, we define the {\bf dual graph}
$\Gamma(P)$ as a formal graph whose set of vertices is $P$ and whose edges
join pairs of pieces that share an arc of external ray. It is always possible
to produce an isomorphic model of $\Gamma(P)$ sitting in the plane, without
intersecting edges and such that it respects the natural immersion of
$\Gamma(P)$ in the plane.

\defn{} When $P = P_n$, we call $\Gamma_n:=\Gamma(P_n)$ the {\bf puzzle graph}
  of depth $n$. In this context, the vertices corresponding to the central
  piece $Y_0^{(n)}$ and the piece around the critical value $f_c(0)$ are
  denoted $\xi_n$ and $\eta_n$ respectively.
\medskip

\begin{center}\begin{figure}[h]
  \includegraphics{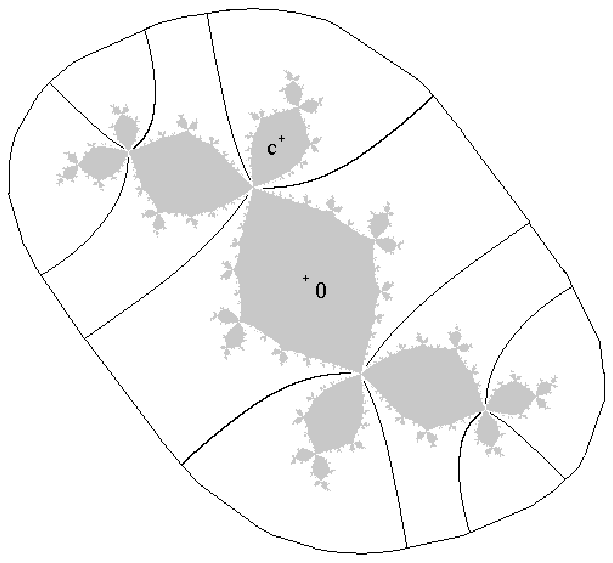}
  \includegraphics{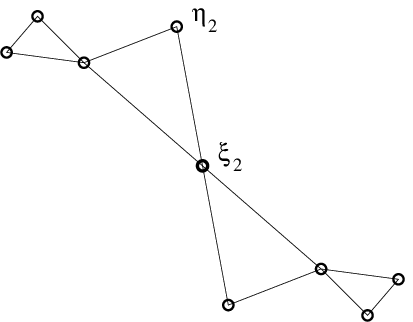}
  \caption[Puzzle of depth 2 and its corresponding graph]
          {\label{fig:Puzzle_Graphs} \it Puzzle of depth 2 and its
          corresponding graph. Splitting the graph at $\xi_2$ we obtain the
          graphs $\text{Puzz}_2^-$ and $\text{Puzz}_2^+$; both shaped like a
          bow tie and isomorphic to $\Gamma_1$. }
\end{figure}\end{center}

\defn{} The vertices $\xi_n$ and $\eta_n$ determine two partial orders on the
  vertex set of $\Gamma_n$ as follows: If $a,b \in V(\Gamma_n)$, we write $a
  \succ_{\eta_n} b$ when every path from $a$ to $\eta_n$ passes through $b$.
  We write $a \succ_{\xi_n} b$ when every path from $a$ to $\xi_n$ passes
  through $b$ or through its symmetric image with respect to the origin.
\medskip

The following are natural consequences of the definitions; see Figure
\ref{fig:Puzzle_Graphs} for reference.

\begin{prop}\label{prop:Puzzle_Properties}
  The puzzle graphs of $f$ satisfy:

  {\renewcommand{\labelenumi}{{\bf G\arabic{enumi}}}
  \renewcommand{\theenumi}{\labelenumi}
  \begin{enumerate}
    \item \label{symmetry} $\Gamma_n$ has 2-fold central symmetry around
          $\xi_n$.

    \item \label{q-gon} $\Gamma_0$ is a $q$-gon whenever $c \in L_{p/q}$. For
          $n \geq 1$, $\Gamma_n$ consists of $2^n$ $q$-gons linked at their
          vertices in a tree-like structure; i.e. the only cycles on this
          graph are the $q$-gons themselves.

    \item \label{next_level} For $n \geq 1$, removing $\xi_n$ and its edges
          splits $\Gamma_n$ into 2 disjoint (possibly disconnected) isomorphic
          graphs. Reattaching $\xi_n$ to each, and adding the corresponding
          edges defines the connected graphs $\text{Puzz}_n^-$ and
          $\text{Puzz}_n^+$ (here, $\eta_n \in \text{Puzz}_n^-$). Then
          $\Gamma_n = \text{Puzz}_n^- \cup \text{Puzz}_n^+$ and
          $\text{Puzz}_n^- ,\, \text{Puzz}_n^+$ are isomorphic to
          $\Gamma_{n-1}$ with $\mp \eta_n$ playing the role of $\xi_{n-1}$ in
          $\text{Puzz}_n^{\pm}$.

    \item \label{maps} For $n \geq 1$ there are two natural maps:
          $f^*:\Gamma_n \longrightarrow \Gamma_{n-1}$ induced by $f$, and
          $\iota^*:\Gamma_n \longrightarrow \Gamma_{n-1}$ induced by the
          inclusion among pieces of consecutive depths. $f^*$ is 2 to 1 except
          at $\xi_n$ and sends $\text{Puzz}_n^{\pm}$ onto $\Gamma_{n-1}$. In
          turn, $\iota^*$ collapses the outermost $q$-gons into vertices.

    \item \label{order} The map $f^*:\big( \Gamma_n,\, \succ_{\xi_n} \big)
          \longrightarrow \big( \Gamma_{n-1},\, \succ_{\eta_{n-1}} \big)$
          respects order. That is, if $a \succ_{\xi_n} b$ then $f^*(a)
          \succ_{\eta_{n-1}} f^*(b)$.
  \end{enumerate}}
\end{prop}

\defn{} Let $\Gamma$ be a graph isomorphic to a subgraph of $\Gamma_n$ and
  $\Gamma'$ a graph isomorphic to a subgraph of $\Gamma_{n-1}$. A map
  $E:\Gamma \longrightarrow \Gamma'$ that satisfies \ref{symmetry} and
  \ref{q-gon} will be called {\bf admissible} if it also respects order in the
  sense of \ref{order}.
\medskip

\pf{ of Proposition \ref{prop:Puzzle_Properties}} Property \ref{symmetry} and
  the existence of $f^*$ and $\iota^*$ are immediate consequences of the
  structure of quadratic Julia sets. The configuration of $\Gamma_0$ is given
  by the rotation number around $\alpha$ and then the tree-like structure of
  $\Gamma_n \, (n \geq 1)$ follows from \ref{next_level}.

  Consider a centrally symmetric simple curve $\gamma \subset Y_0^{(n)}$
  connecting two opposite points of the equipotential curve $e_{(2^{2^{-n}})}$
  that bounds $Y_0^{(n)}$. Then $\gamma$ splits the simply connected region
  $\bigcup_{Y \in P_n} Y$ in 2 identical parts. Therefore, $\Gamma \setminus
  \xi_n$ is formed by 2 disjoint graphs justifying the existence of
  $\text{Puzz}_n^{\pm}$. However, $\bdry Y_0^{(n)}$ may contain several
  segments of $e_{(2^{2^{-n}})}$; so $\gamma$, and consequently
  $\text{Puzz}_n^{\pm}$, {\it are not uniquely determined}. This ambiguity is
  not consequential; Lemmas \ref{lemma:New_Frame_Graph} and
  \ref{lemma:Embedding_of_F} describe the proper method of handling it.

  The fact that $f$ maps the central piece to a non-central one containing the
  critical value legitimizes the selection of $\text{Puzz}_n^-$ as the unique
  graph containing $\eta_n$. By symmetry, every piece of $P_n$ except the
  central one has a symmetric partner and they both map in a 1 to 1 fashion to
  the same piece of $P_{n-1}$. The isomorphisms in \ref{next_level} follow.

  If two pieces $A,B$ of depth $n$ share a boundary ray, their images will
  too. Moreover, letting $A',B'$ be the pieces of depth $n-1$ containing $A$
  and $B$, it is clear that $\bdry A'$ and $\bdry B'$ must share the {\it
  same} ray as $\bdry A$ and $\bdry B$. This shows that $f^*$ and $\iota^*$
  effectively preserve edges and are well defined graph maps. Clearly $f^*$ is
  2 to 1, so to complete the proof of \ref{maps} we only need to justify the
  collapsing property of $\iota^*$, and by Property \ref{next_level}, it is
  sufficient to consider the case $\iota^*:\Gamma_1 \longrightarrow
  \Gamma_0$. Now, the non-critical piece $Y_j^{(0)}$ contains a unique piece
  $Y_j$ of $P_1$. However, the critical piece $Y_0^{(0)}$ contains a total of
  $q$ different pieces of depth 1: a smaller central piece $Y_0^{(1)}$ and
  $q-1$ lateral pieces $-Y_j$. The resulting graph, $\Gamma_1$, consists then
  of two $q$-gons joined at the vertex $\xi_1$. Under $\iota^*$, one of these
  $q$-gons collapses on the critical vertex $\xi_0$.

  To prove \ref{order}, let us construct the tree $\Gamma'_n$ with 2 to 1
  central symmetry by collapsing every $q$-gon into a single vertex. The
  orders $\succ_{\xi'_n},\, \succ_{\eta'_n}$ in $\Gamma'_n$ are induced by the
  orders in $\Gamma_n$. Then the corresponding map ${f^*}':\big( \Gamma'_n,\,
  \succ_{\xi_n} \big) \longrightarrow \big( \Gamma'_{n-1},\, \succ_{\eta_n}
  \big)$ is a 2 to 1 map on trees that takes each half of $\Gamma'_n$
  injectively into a sub-tree of $\Gamma'_{n-1}$ and respects order. Since
  vertices in a cycle are not ordered, $f^*$ respects order as well.
\QED

\subsection{Parapuzzle}
While the puzzle encodes the combinatorial behavior of the critical orbit for
a specific map $f_c$, the {\it parapuzzle} dissects the parameter plane into
regions of parameters that share similar behaviors: In every wake of $M$ we
define a partition in pieces of increasing depths, with the property that all
parameters inside a given {\it parapiece} share the same critical orbit
pattern up to a specific depth.

\defn{} Consider a wake $W_{p/q}$ and let $n \geq 0$ be given. Call $W^n$ the
  wake $W_{p/q}$ truncated by the equipotential $e_{\left( 2^{2^{-n}}
  \right)}$ and consider the set of angles $\mathcal{P}_n(\frac pq) = \big\{ t
  \mid 2^nt \in \mathcal{P}\big( \frac pq \big) \big\}$ (compare Subsection
  \ref{subsect:External_Rays_n_Wakes}). The {\bf parapieces} of $W_{p/q}$ at
  depth $n$ are the closures of the components of $W^n \setminus \big\{r_t
  \mid t \in \mathcal{P}_n(\frac pq) \big\}$.
\medskip

\note{} Even though the critical value $f_c(0)$ is simply $c$, it will be
  convenient to write $c \in \Delta$ when $\Delta$ is a parapiece and $f_c(0)
  \in V$ when $V$ is a piece in the dynamical plane of $f_c$. In general, we
  will use the notation $\text{OBJ}[c]$ to refer to dynamically defined
  objects $\text{OBJ}$ associated to a specific parameter $c$.

\defn{} When the boundary of a dynamical piece $A$ is described by the same
  equipotential and ray angles as those of a parapiece $B$, we denote this
  relation by $\bdry A \circeq \bdry B$.
\medskip

\defn{} Let $c \in M$ be a parameter whose puzzle is defined up to depth $n$.
  We denote by $\text{CV}_n[c] \in P_n[c]$ the piece of depth $n$ that
  contains the critical value: $f_c(0) \in \text{CV}_n[c]$.
\medskip

A consequence of Formula \ref{eqn:Phi_M} is the well known fact that follows.
For a proof of the main statement, refer to \cite{DH_p-l} or
\cite{Roesch}. For a proof of the winding number property, refer to
\cite{Chirurgie} and Proposition 3.3 of \cite{L_parapuzzle}; also, see the
Appendix for the definition of holomorphic motions.

\begin{prop}\label{prop:Identify_Bdries}
  Let $\Delta$ be a parapiece of depth $n$ in some wake $W$. Then
  $\text{CV}_n[c] \circeq \Delta$ for every $c \in \Delta$ so the family
  $\big\{ c \mapsto \text{CV}_n[c] \mid c \in \Delta \big\}$ is well defined;
  it determines a holomorphic motion of the critical value pieces. The
  holomorphic motion has $\big\{ c \mapsto f_c(0) \big\}$ as a section with
  winding number 1.
\end{prop}

We can interpret the result on winding number as loosely saying that, as $c$
goes once around $\bdry \Delta$, the critical value $f_c(0)$ goes once around
$\bdry\text{CV}_n$. However, this description is not entirely accurate since
$\bdry\text{CV}_n[c]$ changes with $c$.

Let us mention the following examples of combinatorial properties that depend
on the behavior of the first $n$ iterates of 0. The fact that these entities
remain unchanged for $c \in \Delta$ follows from Proposition
\ref{prop:Identify_Bdries} and will be useful in the next sections.

\begin{itemize}
  \item The isomorphism type of $\Gamma_n[c]$.
  \item The combinatorial boundary of every piece of depth $\leq n$.
  \item The location within $P_n[c]$ of the first $n$ iterates of the critical
        orbit. \\
\end{itemize}

From the general results of \cite{L_parapuzzle}, we can say more about the
geometric objects associated to the above examples.

\begin{prop}
  Each of the sets listed below moves holomorphically as $c$ varies in
  $\Delta$:

  \begin{itemize}
    \item The boundary of every piece of depth $\leq n$.
    \item The first $n$ iterates of the critical orbit.
    \item The collection of $j$-fold preimages of $\alpha$ and $\beta$ $(j
          \leq n)$.
  \end{itemize}
\end{prop}

\subsection{Principal nest}
The principal nest is well defined for parameters $c$ that belong neither to
$\overline{\cardioid}$ nor to an immediate component. The first condition
means that both fixed points are repelling (so the puzzle is defined), while
the second condition characterizes those polynomials that do not admit an {\it
immediate renormalization} as described below. We restrict further to
parameters $c$ such that the orbit of 0 is recurrent to ensure that the nest
is infinite. These necessary conditions will justify themselves as we describe
the nest. \\

\begin{center}\begin{figure}[h]
  \includegraphics{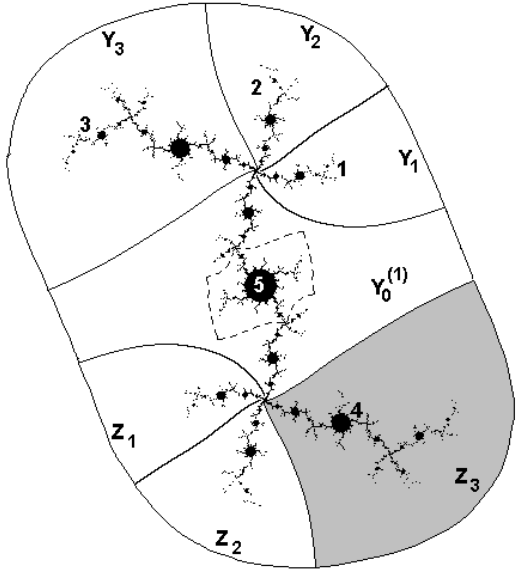}
  \caption[First stage of construction of the principal nest]
          {\label{fig:Nest_Of_Level_0} \it Puzzle $P_1(f_c)$ of depth 1, where
          $c=(0.35926...)+i(0.64251...)$ is the center of the component of
          period 5 in $L_{1/4}$. The first escape is $f_c^{\circ 4}(0) \in
          Z_3$ and the pull back $V_0^0$ is shown in dotted lines. Note that
          $f^{\circ 5}(0) \in V_0^0$. This creates at once the piece $V_0^1
          \Subset V_0^0$ around the central component of $\mathbb{C} \setminus
          J_c$ ($V_0^1$ is not shown).}
\end{figure}\end{center}

In order to explain the construction of the principal nest, we need a more
detailed description of the puzzle partition at depth 1 (use Figure
\ref{fig:Nest_Of_Level_0} for reference). As a note of warning, the pieces of
depth 1 will be renamed to reflect certain properties of $P_1$. That is, we
will override the use of the symbols $Y_j^{(1)}$.

The puzzle depth $P_1$ consists of $2q-1$ pieces of which $q-1$ are the
restriction to lower equipotential of the pieces $Y_1^{(0)}, Y_2^{(0)},
\ldots, Y_{q-1}^{(0)}$. Such pieces cluster around $\alpha$ and will be
denoted $Y_1, Y_2, \ldots, Y_{q-1}$. The restriction of $Y_0^{(0)}$ however,
is further divided into the union of the critical piece $Y_0^{(1)}$ and $q-1$
pieces $Z_1, Z_2, \ldots, Z_{q-1}$ which are symmetric to the corresponding
$Y_j$ and cluster around $-\alpha$. The indices are again determined by the
rotation number of $\alpha$ so that $f(Z_j)$ is opposite to $Y_j$ and
consequently $f(Z_j) = Y_{j+1}^{(0)}$.

Note that $f^{\circ q}(0) \in Y_0^{(0)}$, so we face two possibilities. It may
happen that $f^{\circ jq}(0) \in Y_0^{(1)}$ for all $j$, in which case we can
find {\it thickenings} of $Y_0^{(1)}$ and $Y_0^{(0)}$, that yield the {\bf
immediate renormalization} $f^{\circ q}: Y_0^{(1)} \longrightarrow Y_0^{(0)}$
described by Douady and Hubbard; or else, we can find the least $k$ for which
the orbit of $0$ under $f^{\circ q}$ escapes from $Y_0^{(1)}$. We will assume
that this is the case, so $f^{\circ kq}(0) \in Z_\nu$ for some $\nu$ and we
call $kq$ the {\bf first escape time.}

The initial nest piece $V_0^0$ is defined as the $(kq)$-fold pull back of
$Z_{\nu}$ along the critical orbit; that is, the unique piece that satisfies
$0 \in V_0^0$ and $f^{\circ kq}(V_0^0) = Z_{\nu}$. In fact, $V_0^0$ can also
be defined as {\bf the largest central piece that is compactly contained in}
$Y_0^{(1)}$: Notice that $Z_\nu \Subset Y_0^{(0)}$ so $V_0^0 \Subset
Y_0^{(1)}$; that is, $\big( \text{int}\, Y_0^{(1)} \big) \setminus V_0^0$ is a
non-degenerate annulus.

The higher levels of the principal nest are defined inductively. Suppose that
the pieces $V_0^0, V_0^1, \ldots, V_0^n$ have been already constructed. If the
critical orbit never returns to $V_0^n$ then the nest is finite. Otherwise,
there is a first return time $\ell_n$ such that $f^{\circ \ell_n}(0) \in
V_0^n$; then we define $V_0^{n+1}$ as the {\it critical} piece that maps to
$V_0^n$ under $f^{\circ \ell_n}$.

\begin{prop} The principal nest $V_0^0 \Supset V_0^1 \Supset \ldots$ is a
  family of strictly nested pieces centered around 0.
\end{prop}

\pf{} $V_0^0$ is a piece of depth $kq$ (the first escape time). Since $V_0^1$
  is a $f^{\circ \ell_1}$-pull back of $V_0^0$, it is a piece of depth
  $kq+\ell_1$ and, in general, $V_0^n$ will be a piece of depth
  $kq+\ell_1+\dots+\ell_n$. Since all pieces contain 0, Property
  \ref{piece_relns} implies that $V_0^j \supset V_0^{j+1}$.

  Recall that $V_0^0 \Subset Y_0^{(1)}$; thus, the $f^{\circ \ell_1}$-pull
  backs of these 2 pieces satisfy $V_0^1 \Subset X$ with $X$ a central piece
  of depth $1+\ell_1$. Now, $0 \nin Z_{\nu}$, so $f^{\circ kq}(0)$ requires
  further iteration to reach a central piece; i.e., $\ell_1 > kq$. By
  construction, $V_0^0$ is a central piece of depth $1+kq$, so Property
  \ref{piece_relns} implies $V_0^1 \Subset X \subset V_0^0$. An analogous
  argument yields the strict nesting property for the nest pieces of higher
  depth.
\QED

\defn{} The {\bf principal annuli} $V_0^{n-1} \setminus V_0^n$ will be denoted
  $A_n$.
\medskip

It may happen that $\ell_{n+1} = \ell_n$; this means that not only does 0
return to $V_0^n$ under $f^{\circ \ell_n}$, but even deeper to $V_0^{n+1}$
without further iteration. In this case we say that the return is {\bf
central} and we call a chain of consecutive central returns $\ell_n =
\ell_{n+1} = \ldots = \ell_{n+s}$ a {\bf cascade of central returns}. An
infinite cascade means that the sequence $\{ \ell_n \}$ is eventually
constant, so $f^{\circ \ell_n}(0) \in \bigcap_{j=n}^\infty V_0^j$. By
definition, $f^{\circ \ell_n}:V_0^{n+1} \longrightarrow V_0^n$ is a {\bf
renormalization} of $f$; that is, a 2 to 1 branched cover of $V_0^n$ such that
the orbit of the critical point is defined for all iterates.

The return to $V_0^n$, however, can be non-central. In fact, it is possible to
have several returns to $V_0^n$ before the critical orbit hits $V_0^{n+1}$ for
the first time. When a return is non-central, the description of the nest at
that level is completed by the introduction of the {\bf lateral} pieces $V_k^n
\in V_0^{n-1} \setminus V_0^n$. Let $\mathcal{O} \subset K$ denote the
critical orbit $\mathcal{O} = \left\{ f^{\circ j}(0)|j \geq 0 \right\}$ and
take a point $z \in \overline{\mathcal{O}} \cap V_0^{n-1}$ whose forward orbit
returns to $V_0^{n-1}$. If we call $r_{n-1}(z)$ the first return time of $z$
back to $V_0^{n-1}$, we can define $V^n(z)$ as the unique puzzle piece that
satisfies $z \in V^n(z)$ and $f^{\circ r_{n-1}(z)}\big(V^n(z)\big) =
V_0^{n-1}$. In particular, it is clear that $V^n(0)$ is just the same as
$V_0^n$ and that any 2 pieces created by this process are disjoint or
equal.

\defn{} The collection of all pieces $V^n(z)$ for $z \in
  \overline{\mathcal{O}} \cap V_0^{n-1}$ that actually contain a point of
  $\mathcal{O}$ is denoted $\mathcal{V}^n$ and referred to as the {\bf level}
  $n$ of the nest.
\medskip

\begin{center}\begin{figure}[h]
  \includegraphics{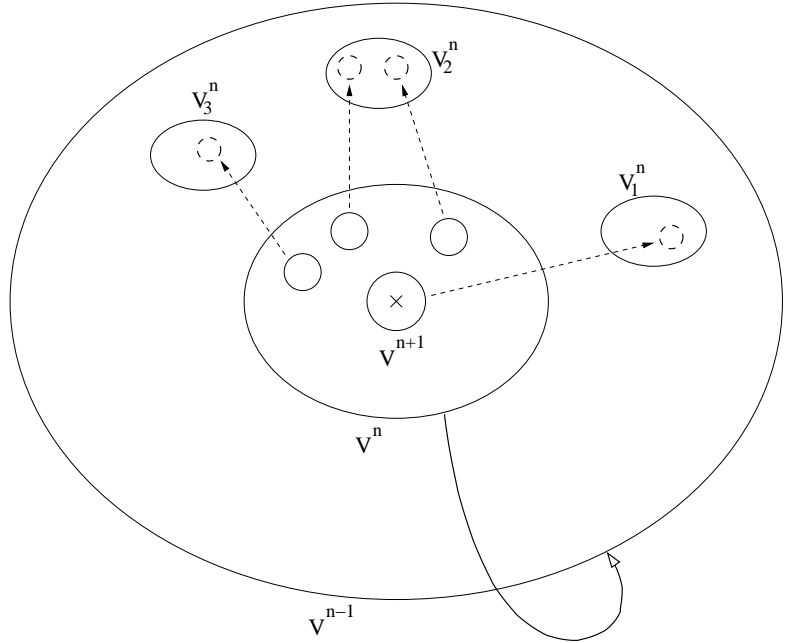}
  \caption[Relation between consecutive nest levels]{\label{fig:Sample_Nest}
           \it Relation between consecutive nest levels. The curved arrow
           represents the first return map $f^{\circ \ell_n}:V_0^n
           \longrightarrow V_0^{n-1}$ which is 2 to 1. The dotted arrows show
           a possible effect of this map on each nest piece of level $n+1$.
           Each $V_j^{n+1}$ may require a different number of additional
           iterates to return to this level and map onto $V_0^n$.}
\end{figure}\end{center}

Under the assumption that $c$ is recurrent, the principal nest will have
infinitely many levels. Let us assume the parameter $c$ is not periodic. Then
it is called {\bf reluctantly recurrent} if for some central piece $V_0^n$
there are arbitrarily long sequences of univalent $f_c$-pull backs of $V_0^n$
along backward orbits in the postcritical set
$\overline{\mathcal{O}}$. Otherwise, $c$ is called {\bf persistently
recurrent}.

\begin{lemma}{\rm {(see \cite{L_attractor},\cite{Martens})}}
\label{lemma:Persistent_Reluctant}
  If $f_c$ is persistently recurrent, $\overline{\mathcal{O}}$ is a Cantor set
  and the action of $f_{c|_{\overline{\mathcal{O}}}}$ is minimal. When $f_c$
  is not renormalizable, $c$ is reluctantly recurrent if and only if some
  central piece $V_0^n$ has infinitely many 1 to 2 pull backs along backward
  orbits of $\mathcal{O}$.
\end{lemma}

\obs{} In particular, if $c$ is non-renormalizable but every level of the
  principal nest has a finite number of pieces, then $f_c$ acts minimally on
  the postcritical set. In this situation, we can name the pieces
  $\mathcal{V}^n = \{ V_0^n,V_1^n, \ldots ,V_{m_n}^n \}$ in such a way that
  the first visit of the critical orbit to $V_i^n$ occurs before the first
  visit to $V_j^n$ whenever $i<j$. Obviously, the value of $r_{n-1}(z)$ is
  independent of $z \in V_k^n$; thus we will denote it $r_{n,k}$.

\defn{}\label{defn:GQL_Maps} For finite $\mathcal{V}^n$ we define the map:
  $$
    g_n:\bigcup_{\mathcal{V}^n} V_k^n \longrightarrow V_0^{n-1},
  $$
  given on each $V_k^n$ by ${g_n}_{|_{V_k^n}} \equiv f^{\circ r_{n,k}}$.
\medskip

The map $g_n$ satisfies the properties of a {\it generalized quadratic-like}
({\bf gql}) {\it map,} i.e.:

\begin{itemize}
  \item $|\mathcal{V}^n| < \infty$.

  \item $\bigcup_{\mathcal{V}^n} V_k^n \Subset V_0^{n-1}$ and all the pieces
        of $\mathcal{V}^n$ are pairwise disjoint.

  \item ${g_n}_{|_{V_k^n}} : V_k^n \longrightarrow V_0^{n-1}$ is a 2 to 1
        branched cover or a conformal homeomorphism depending on whether $k=0$
        or not.
\end{itemize}

Note that $g_n$ usually is the result of a different number of iterates of $f$
when restricted to different $V_k^n$. However, since we often refer to the map
$g_n$ as acting on individual pieces, it is typographically convenient to
introduce the notation

\defn{}\label{defn:g_{n,k}} The map ${g_n}_{|_{V_k^n}} = f^{\circ r_{n,k}}$
  will be denoted $g_{n,k}$.
\medskip

Thus, $g_{n,k}(V_k^n) = V_0^{n-1}$ is a 2 to 1 branched cover or a
homeomorphism depending on whether $k=0$ or not.

From this moment on, we will assume that the principal nest is infinite, and
that $f$ is non-renormalizable, thus excluding the possibility of an infinite
cascade of central returns. In this situation we say that $f$ is {\bf
combinatorially recurrent.}

\subsection{Paranest}
The {\it paranest} is well defined around parameters $c$ outside the main
cardioid that are neither immediately renormalizable nor postcritically
finite.

\defn{}\label{defn:Paranest} If $c$ is a parameter such that $f_c$ has a well
  defined nest up to level $n$ (for $n \geq 0$), the {\bf paranest} piece
  $\Delta^n[c]$ is defined by the condition $\bdry \Delta^n[c] \circeq \bdry
  f_c(V_0^n)$; where $V_0^n$ is the central piece of level $n$ in the
  principal nest of $f_c$. By the Douady-Hubbard theory, $\Delta^n[c]$ is a
  well defined region.
\medskip

The definition of principal nest, together with Proposition
\ref{prop:Identify_Bdries} imply that when $c' \in \Delta^n[c]$, the principal
nests of $f_c$ and $f_{c'}$ are identical until the first return $g_n(0)$ to
$V_0^{n-1}$ (which creates $V_0^n$). In fact, the relevant pieces move
holomorphically as $c'$ varies and $\Delta^n[c]$ is the largest parameter
region over which the initial set of $\ell_n$ iterates of 0 (recall that $g_n
\equiv f^{\circ \ell_n}$) moves holomorphically without crossing piece
boundaries.

Following the presentation of \cite{L_parapuzzle}, the family
$\big\{g_n[c']:V_0^n[c'] \longrightarrow V_0^{n-1}[c'] \mid c' \in \Delta^n[c]
\big\}$ is a proper DH quadratic-like family with winding number 1. The last
property follows from Proposition \ref{prop:Identify_Bdries} since $g_n$ is
the first return to a critical piece at this level. \\

Since the central nest pieces are strictly nested, the above definition
implies that the pieces of the paranest are strictly nested as well. It
follows that $\big( \text{int}\, \Delta^n \big) \setminus \Delta^{n-1}$ is a
non-degenerate annulus. One of the main concerns is to estimate its modulus
or, as it is sometimes called, the {\bf paramodulus.}

\section{Frame system}\label{sect:Frames}
\setcounter{equation}{0}

Let $f_c$ have an infinite principal nest. For real parameters, Lyubich
provides in \cite{L_attractor} a complete criterion for compatibility between
consecutive nest levels. Since the Julia set is an interval when $c \in
\mathbb{R}$, the compatibility conditions are given in terms of the left/right
location of lateral pieces (relative to 0) and the orientation of each
$g_{n,k}$ (as an interval map).

In the case of a complex parameter, the nest falls short of being a complete
invariant for the dynamics of the critical orbit. The reason is that the nest
description does not account for the relative positions between lateral
pieces. In contrast to the real case, the Julia set of a complex polynomial
displays a complicated structure that varies with the parameter. Lateral
pieces may be attached to different branches of the Julia set. For this
reason, a record of the relative positions of nest pieces must be preceded by
a description of the combinatorial structure around them.

In this Section we enhance the principal nest with the addition of a {\it
frame system.} This provides the necessary language to locate the lateral nest
pieces and describe as a consequence, the behavior of the critical orbit. {\it
The idea is to split the central nest pieces in smaller regions by a procedure
that resembles the construction of the puzzle.} \\

For convenience, let us summarize certain aspects of the construction before
giving it in detail. Recall that the definition of $V_0^0$ guarantees that
$\big( \text{int}\, Y_0^{(1)} \big) \setminus V_0^0$ is a non-degenerate
annulus. Because of this initial step, and since our purpose is that frame
levels correspond to nest levels, we need to pay individual attention to the
construction of the first three levels of the frame. Figure
\ref{fig:Two_Frames} illustrates these initial steps. We will keep in mind our
convention of distinguishing between puzzle depths and nest
levels. Accordingly, frames will be also stratified in levels since their
definition depends on the same pull backs as those used for the nest. To
distinguish between nest pieces and frame pieces, the latter will be referred
to as {\it cells}.  As a final note of warning, we will abuse our notation and
use $F_n$ to refer to the frame as well as to the system of curves that bound
its cells. In particular, we will use $\bdry F_n$ to describe the union of
curves that form the boundary of the union of all cells in $F_n$. The context
will always make clear which meaning is intended.

\subsection{Frames}
As mentioned above, some attention must be given to the construction of the
frames $F_0, F_1$ and $F_2$ so that the properties in Proposition
\ref{prop:Frame_Properties} hold. Figure \ref{fig:Nest_Of_Level_0} provides a
useful reference. After this, the frames of higher levels are defined
inductively.

Consider the puzzle partition at depth 1 and recall that $kq$ denotes the
first escape of the critical orbit to $Z_{\nu}$. The {\bf initial frame} $F_0$
is the collection of nest pieces $F_0 = \big\{Y_0^{(1)} \big\} \cup \big\{
\bigcup_{j=1}^q \{Z_j \} \big\}$, each of which is called a {\bf frame
cell}. In particular, $\Gamma(F_0)$ is a $q$-gon. The frame $F_1$ is the
collection of $f^{\circ kq}$-pull backs of cells in $F_0$ along the orbit of
0.

From the definition, one of the cells of $F_1$ is the central piece $V_0^0$
that maps 2 to 1 onto $Z_{\nu} \in F_0$. The pull back of any other cell $A
\in F_0$ consists of two symmetrically opposite cells, each mapping
univalently onto $A$. We say that $F_1$ is a {\it well defined unimodal} pull
back of $F_0$.

\begin{lemma}\label{lemma:Start_Frame}
  All the cells of $F_1$ are contained in $Y_0^{(1)}$.
\end{lemma}

\pf{} Since $kq > 1$, $f^{\circ kq}\big( Y_0^{(1)} \big)$ is an extension of
  $Y_0^{(0)}$ to a larger equipotential. Thus, $f^{\circ kq}\big( Y_0^{(1)}
  \big)$ contains all cells of $F_0$.
\QED

Let $\lambda$ be the first return time of 0 to a cell of $F_1$. By Lemma
\ref{lemma:Start_Frame}, the collection $F'_2$ of pull backs of cells in
$F_1$ along the $f^{\circ \lambda}$-orbit of 0 is well defined and 2 to 1.
Unfortunately, it does not cover every point of $J_f$ inside $V_0^0$. We will
give first some results about $F'_2$ and define afterward a complete frame
of level 2.

\begin{lemma}
  The temporary frame $F'_2$ satisfies:
  \begin{enumerate}
    \item All cells of $F'_2$ are contained in $V_0^0$.
    \item $V_0^1$ is contained in the central cell of $F'_2$.
  \end{enumerate}
\end{lemma}

\pf{}
  First note that $\lambda = kq + (q - \nu)$ is the first return of 0 to
  $Y_0^{(1)}$ after the first escape to $Z_{\nu}$. We have
  $kq < \lambda \leq \ell_0$, where the second inequality is true since
  $V_0^0 \in F_1$. Then the first return to $F_1$ occurs no later than the
  first return to $V_0^0$. By definition, $f^{\circ \lambda}(V_0^0)$ is
  just $Y_0^{(0)}$ extended to a larger equipotential. Since all cells of
  $F_1$ are inside $Y_0^{(1)} \subset f^{\circ \lambda}(V_0^0)$, the first
  assertion follows.

  Now, $V_0^1$ is central. By the Markov properties of $\mathcal{Y}_c$,
  either $V_0^1$ is contained in the central cell $C$ of $F'_2$ or vice versa.
  However, both $f^{\circ \ell_0}(V_0^1)$ and $f^{\circ \lambda}(C)$
  belong to $F_1$. Since $\ell_0 \geq \lambda$, the first possibility
  is the one that holds. This proves property (2).
\QED

Our intention is to extend $F'_2$ to a frame that covers the intersection
$J_f \cap V_0^0$. To do this, we just need to add the
$f^{\circ \lambda}$-pull backs of the pieces $Z_{\nu}$. The union of those
pull backs with the cells of $F'_2$ is the frame $F_2$. 

\begin{figure}[h]
\begin{center}
  \includegraphics{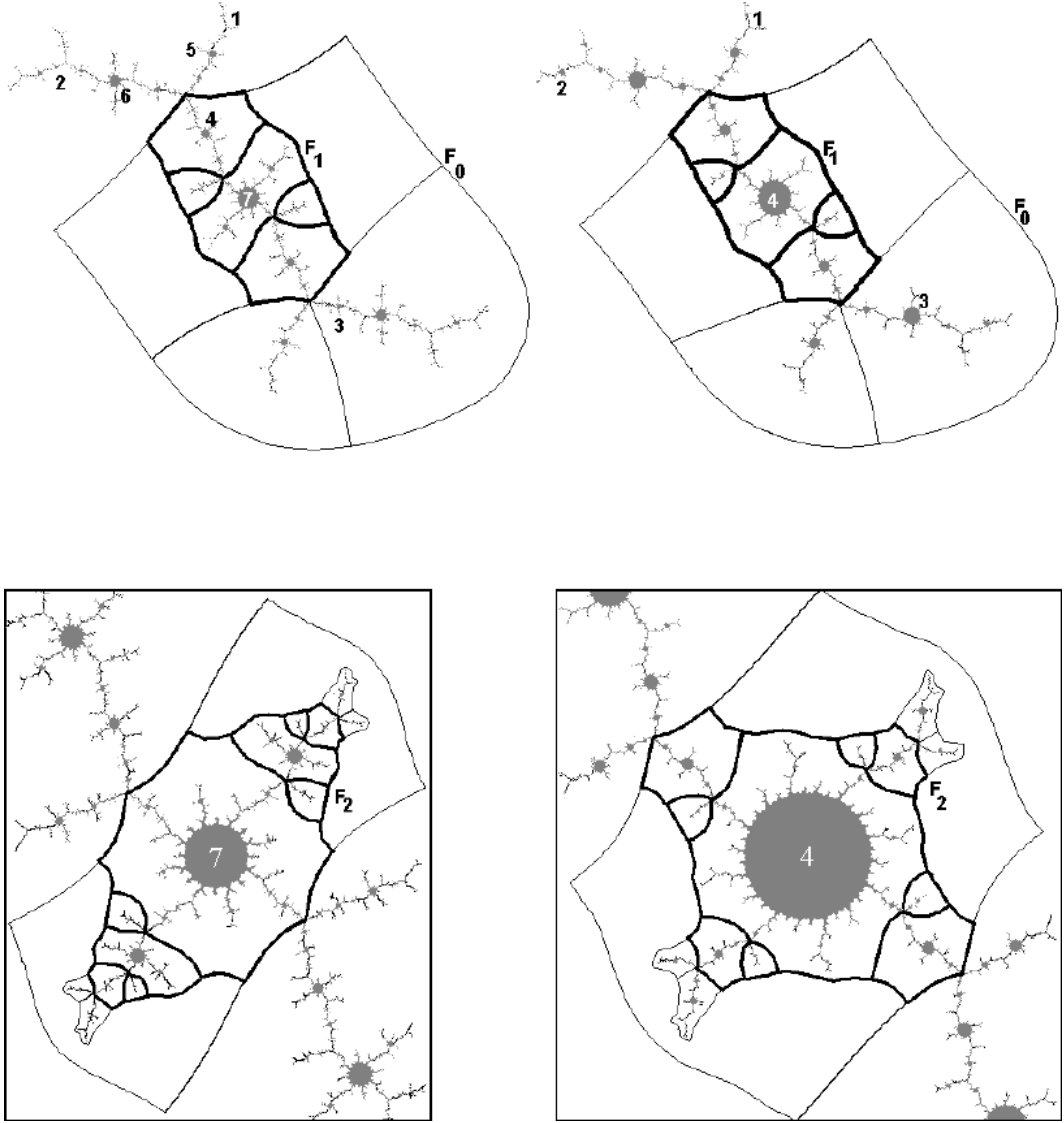}
  \caption[Comparison between two frames with similar combinatorics]
          {
           \label{fig:Two_Frames} \it Both of these parameters belong to the
           left antenna of $L_{1/3}$; they are centers of components of
           periods 7 and 4. Above we can see that the structures of the
           frames of levels 0 and 1 coincide between the two examples. Still,
           the first return to $F_1$ falls in each case on a different cell,
           producing dissimilar frames of level 2. The pull back of cells in
           $F_1$ produces a preliminary frame $F'_2$, shown in heavy line on
           the second row. The complete frame $F_2$, inside $V_0^0$ has
           $2(q-1)$ additional cells (here $q=3$) in order to cover all of
           $J_f \cap V_0^0$.
          }
\end{center}
\end{figure}

After introducing the first frames and relating them to the initial levels
of the nest, we can give the complete definition of the {\it frame system.}
The driving idea of this discussion is that the internal structure of a
frame $F_{n+2}$, represented by the graph $\Gamma(F_{n+2})$, provides a
decomposition of $J_f \cap V_0^n$ that helps to describe the combinatorial
type of the nest at level $n+1$.

\defn{} For $n \geq 0$ consider the first return $g_n(0) \in V_0^n$ and define
  $F_{n+3}$ as the collection of $g_n$-pull backs of cells in $F_{n+2}$ along
  the critical orbit. The family $\mathcal{F}_c = \{F_0,F_1, \ldots \}$ is
  called a {\bf frame system} for the principal nest of $f_c$ and each piece
  of a frame is called a {\bf cell}.
\medskip

  The dual graph $\Gamma(F_n)$ (see Subsection \ref{subsect:Graphs}) is called
  the {\bf frame graph.} As in the case of the puzzle graph, we consider
  $\Gamma(F_n)$ with its natural embedding in the plane.

Let us mention now some properties of frame systems.

\begin{prop} \label{prop:Frame_Properties} The frame system satisfies:
  \begin{enumerate}
    \item Frames exist at all levels.
    \item The union of cells $\bigcup_{C_i \in \mathcal{F}_n} C_i$ forms a
          cover of $K_{f_c} \cap V_0^{n-2}$.
    \item The central cell of $F_n$ contains the nest piece $V_0^{n-1}$.
    \item Each $F_n$ has 2-fold central symmetry around 0.
    \item \label{small_inside} Suppose there is a non-central return; then,
          eventually all nest pieces are compactly contained in cells of the
          corresponding frame.
  \end{enumerate}
\end{prop}

The following observation will help clarify the definition of frames (also,
refer to Figure \ref{fig:Two_Frames}). As follows from the comment after Lemma
\ref{lemma:Start_Frame}, the union of cells in $\mathcal{F}_2$ covers exactly
the intersection of $K_f$ with the nest piece $V_0^0$. This is because $V_0^0$
can be described as the pull back of $Y_0^{(0)}$ under the first return map to
$F_1$. Then, we can think of this union of cells as a single piece, determined
by the same rays as $V_0^0$, but cut off by a lower equipotential. \\

\pf{ of Proposition \ref{prop:Frame_Properties}} $F_0$ and $F_1$ are easily
  seen to exist from their construction. Since $F_1$ covers the central part
  of $K_f$ between $\alpha$ and $-\alpha$, there will definitely be a return
  to it, creating $F'_2$. As we saw already, this frame is contained inside
  $V_0^0$, so its pull backs are well defined as long as there are new levels
  of the nest. In particular, this already proves claim 2. Since the principal
  nest is infinite, the critical point is recurrent or the map is
  renormalizable. Either case creates critical returns to central nest pieces
  of arbitrarily high level, so $F_{n+1}$ is defined.

  The piece $V_0^0$ is actually the central cell of $F_1$. Now, the first
  return to $F_1$ cannot occur later than the first return to $V_0^0$, so the
  central cell $C$ of $F_2$ is of lower depth than $V_0^1$; thus, $V_0^1
  \subset C$. Afterwards, the depth from $V_0^{n-1}$ to $V_0^n$ increases by
  $\ell_{n-1}$, while the depth from $F_n$ to $F_{n+1}$ increases
  $\ell_{n-2}$. Inductively, since $V_0^{n-1} \subset F_n$ and $\ell_{n-2}
  \leq \ell_{n-1}$, we obtain $V_0^n \subset F_{n+1}$.

  Now, each $F_n$ is a well defined 2 to 1 pull back of $F_{n-1}$, so a cell
  $C$ belongs to $F_n$ if and only if its symmetric $-C \in F_n$. Finally,
  Part (\ref{small_inside}) follows in a similar manner to the analogous
  property of $V_0^0$ inside $Y_0^{(0)}$.
\QED

\subsection{Frame labels}
Our next objective is to introduce a labeling system for pieces of the
frame. This will allow us to describe the relative position of pieces of the
nest within a central piece of the previous level.  Unlike the case of
unimodal maps, where nest pieces are always located left or right of the
critical point, the possible labels for vertices of $\Gamma(F_n)$ will depend
on the combinatorics of the critical orbit. Only after determining the
labeling, it becomes possible to describe the location of nest pieces in a
systematic manner.

Observe that the structure of $F_{n+1}$ is trivially determined once we know
$F_n$ and the location of $g_n(0)$. A graphic way of seeing this is as
follows. Say that the first return $g_{n-1}(0)$ to $V_0^{n-2}$ falls in a cell
$X \in F_n$. Let $L_n$ and $R_n$ be two copies of $\Gamma(F_n)$ with disjoint
embeddings in the plane. Now connect $L_n$ and $R_n$ with a curve $\gamma$
that does not intersect either graph. Suppose that one extreme of $\gamma$
lands at the vertex of $L_n$ that corresponds to $X$ and the other extreme
lands at the corresponding vertex of $R_n$ {\it approaching it from the same
access}.

\begin{figure}[h]
\begin{center}
  \includegraphics{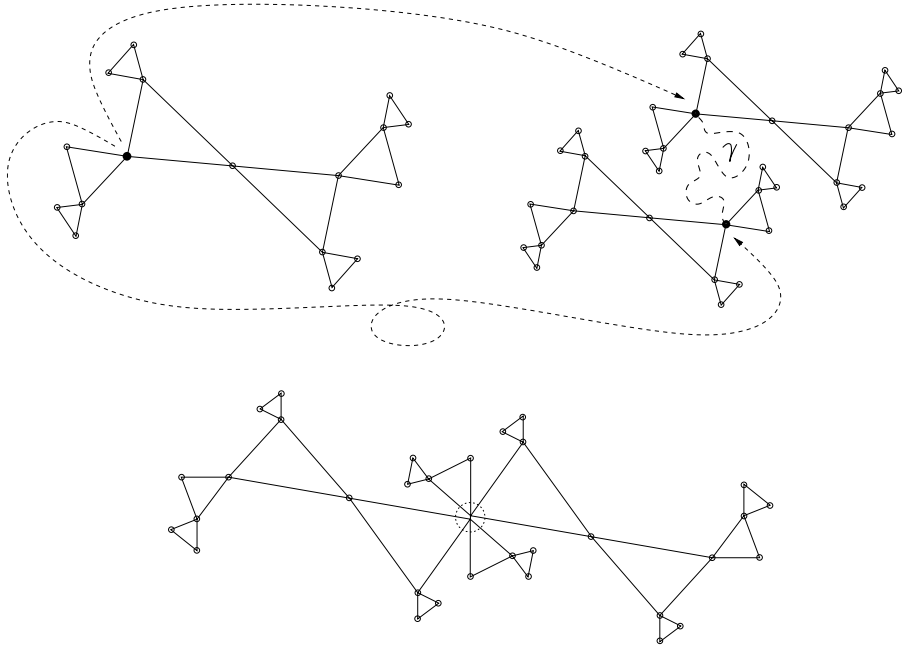}
  \caption[Constructing the next level graph.]{\label{fig:New_Frame} \it The
           curve $\gamma$ joins two copies of the same frame graph
           approaching the selected vertex from the same direction. The new
           frame graph is obtained after $\gamma$ is contracted to a point.}
\end{center}
\end{figure}

\begin{lemma}\label{lemma:New_Frame_Graph}
  If $\gamma$ is collapsed by a homotopy of the whole ensemble, the resulting
  graph is isomorphic to $\Gamma(F_{n+1})$.
\end{lemma}

\note{} The above construction provides $\Gamma(F_{n+1})$ with a natural plane
  embedding; see lemma \ref{lemma:Embedding_of_F} below. \\

A label at level $n$ will be a chain of $n+1$ symbols taken from the alphabet
\{ {\sf 0,1, \ldots, (q-1), l, r, e, b, t} \}. First, put the labels
\{ {\sf '0', '1', \ldots , '(q-1)'} \} on the cells of $F_0$, starting at the
central piece $Y_0^{(0)}$ and moving counterclockwise.

Let $\sigma_0$ be the label of the cell that holds the first return of 0 to
$F_0$ and, in general, let $\sigma_n$ denote the label of the cell in
$\Gamma(F_n)$ that holds the first return of 0. In order to label
$\Gamma(F_{n+1})$, assume that we know the number $q$ of pieces in $F_0$, and
the {\it label sequence} $(q;\sigma_0, \ldots , \sigma_{n-1})$ that identify
the location of first returns of 0 to levels $0, \ldots , n-1$ of the nest.
In particular, all frames up to $\Gamma(F_n)$ have been successfully labeled.

Duplicate in $L_n$ the labels of $\Gamma(F_n)$, but concatenate an extra
{\sf 'l'} at the beginning. Do a similar labeling on $R_n$ by concatenating an
extra {\sf 'r'} to the duplicated labels. Note that the labels of the two
vertices corresponding to $X$ are ${\sf 'l'}\sigma_n$ and ${\sf 'r'}\sigma_n$.
The labels on $\Gamma(F_{n+1})$ will be the same as those in the union of $L_n$
and $R_n$ except that we change the label of the identified vertex, to become
${\sf '0'}\sigma_n$.

\note{} The above procedure does not give labels to the additional cells of
  $F_2$ that do not come from a pull back. These are the cells that are not
  drawn in heavy line in Figure \ref{fig:Two_Frames} Being cells of level 2,
  their labels should have 3 symbols for consistency with the rest. The
  easiest way to do this is simply to impose the labels {\sf 'et1', 'et2',
  \ldots, 'et(q-1)'} and {\sf 'eb1' ,'eb2', \ldots, 'eb(q-1)'} in their
  natural order in the plane ({\sf 'et'} stands for {\bf e}xtra piece on {\bf
  t}op and {\sf 'eb'} for {\bf e}xtra piece on {\bf b}ottom), then extend the
  labeling to higher levels as described. \\

Clearly, $f$ induces a map $f_*:\Gamma(F_{n+1}) \longrightarrow \Gamma_n$ for
$n \geq 2$, that acts by forgetting the leftmost symbol of each label. This
is the case also for the induced map on the temporary frame $F'_2$.

\subsection{Properties of frame labellings}
Under certain conditions, label sequences give a complete characterization of
the entire combinatorial structure. This is the content of Theorem
\ref{thm:Frame_Renormalizations}. Before stating it, we need to review some
properties of the frame and its labels.

\begin{lemma}\label{lemma:Embedding_of_F}
  The plane embedding of $\Gamma$ does not depend on the homotopy class of
  the curve $\gamma$ in lemma \ref{lemma:New_Frame_Graph}.
\end{lemma}

\pf{} Since we regard $\Gamma = \Gamma(F_n)$ as embedded in the sphere, the
  exterior of $\Gamma$ is simply connected, so there is a natural cyclic order
  of accesses to vertices (some vertices can be accessed from more than one
  direction). In this order, all accesses to $L_n$ are grouped together,
  followed by the accesses to $R_n$.
\QED

It is important to mention that the resulting labeling of $\Gamma(F_n)$ {\bf
does} depend on the access to $\xi_n$ approached by $\gamma$. However, the
final unlabeled graphs are equivalent as embedded in the plane.

As we just mentioned, some vertices are accessible from $\infty$ in two or
more directions. These are precisely the vertices whose label contains the
symbol {\sf '0'} (for $n \geq 1$). Since such a vertex represents a frame cell
that maps (eventually) to a central frame cell, the tail of a label with {\sf
'0'} at position $j$ must be $\sigma_j$. On the other hand, for every $j$
there must be labels with a {\sf '0'} in position $j$. It follows that the set
of labels of $\Gamma(F_n)$ and the sequence $(q;\sigma_0, \ldots , \sigma_n)$
can be recovered from each other.

\subsection{Frames and nest together}
The definition of frame system was conceived to satisfy the properties of
Proposition \ref{prop:Frame_Properties}. An extension of the argument used to
prove those properties shows that every piece $V_j^n$ of the nest is contained
in a frame cell of level $n+1$. Moreover, we would like to extend the
definition of frames so that each $V_j^n$ can be partitioned by a pull back of
an adequate central frame. For this, we must recall first that $g_{n,j}(V_j^n)
= V_0^{n-1} \supset F_{n+1}$.

\defn{} The frame $F_{n,k}$ is the collection of pieces inside $V_k^{n-2}$
  obtained by the $g_{n-2,k}$-pull back of $F_{n-1}$. Elements of the frame
  $F_{n,k}$ are called {\bf cells} and we will write $F_{n,0}$ instead of
  $F_n$, when there is a need to stress that a property holds in $F_{n,k}$ for
  every $k$.
\medskip

  If a puzzle piece $A$ is contained in a cell $B \in F_{n,k}$, we denote $B$
  by $\Phi_{n,k}(A)$.

We have described already how to label $F_n$. The other frames $F_{n,k} \, (k
\geq 1)$, mapping univalently onto $F_{n-1}$, have a natural labeling induced
from that of $F_{n-1}$ by the corresponding $g_{n-2,k}$-pull back.

Let us describe now the itinerary of a piece $V_j^n$. Since $V_j^n \subset
V_0^{n-1}$, the map $g_{n-1}$ takes $V_j^n$ inside some piece
$V_{k_1(j)}^{n-1} \subset V_0^{n-2}$. Then, $g_{n-1,k_1(j)}$ takes
$g_{n-1}(V_j^n)$ inside a new piece $V_{k_2(j)}^{n-1}$ and so on, until the
composition of returns of level $n-1$
$$
  (g_{n-1,k_r(j)} \circ \ldots \circ g_{n-1,k_1(j)} \circ g_{n-1})_{|_{V_j^n}}
$$
is exactly $g_{n,j}:V_j^n \mapsto V_0^{n-1}$. Of course, $k_r$ is just 0, and
we will write it accordingly.

We have extra information that deems this description more accurate. For the
sake of typographical clarity, we will write $k_i$ instead of $k_i(j)$. For $i
\leq r$, let $\Phi_{n+1,k_i}$ be the cell in $F_{n+1,k_i} \subset
V_{k_i}^{n-1}$ that contains
$$
  g_{n-1,k_i} \circ \ldots \circ g_{n-1,k_1} \circ g_{n-1}(V_j^n)
$$
and denote by $\lambda_{n+1,k_i}$ the label of $\Phi_{n+1,k_i}$.

\defn{} The {\bf itinerary} of $V_j^n$ is the list of piece-label pairs:
  \begin{equation}
    \chi(V_j^n) = \left(
    [V_{k_1}^{n-1}; \lambda_{n+1,k_1}],
    [V_{k_2}^{n-1}; \lambda_{n+1,k_2}], \ldots,
    [V_{k_{r-1}}^{n-1}; \lambda_{n+1,k_{r-1}}],
    [V_0^{n-1}; \lambda_{n+1,0}] \right)
  \end{equation}
  up to the moment when $V_j^n$ maps onto $V_0^{n-1}$.
\medskip

Note first of all that the last label, $\lambda_{n+1,0}$, will start with {\sf
'0'} due to the fact that $V_0^{n-1}$ is in the central cell of $F_n$. More
importantly, the conditions
\begin{equation}\label{eqn:Admissibility}
  \begin{array}{lcll}
    V_{k_1}^{n-1} & \subset & g_{n-1}(\Phi_{n+1,0}) & \\
    V_{k_{i+1}}^{n-1} & \subset & g_{n-1,k_i}(\Phi_{n+1,k_i}) &
    2 \leq i <r
  \end{array}
\end{equation}
must hold since we know that
$g_{n-1,k_{i-1}} \circ \ldots \circ g_{n-1,k_1} \circ g_{n-1}(V_j^n)
 \subset \Phi_{n+1,k_i}$ and
$g_{n-1,k_i} \circ \ldots \circ g_{n-1,k_1} \circ g_{n-1}(V_j^n)
 \subset V_{k_{i+1}}^{n-1}$.

\defn{} When we specify the sequence of frame labellings up to a given level
  $n$, the locations of the nest pieces and their (admissible) itineraries, we
  say that we have described the {\bf combinatorial type} of the map at
  level $n$. If $| \mathcal{V}^n | < \infty$ we say that the type is {\bf
  finite}; refer to Lemma \ref{lemma:Persistent_Reluctant} and Definition
  \ref{defn:GQL_Maps}.
\medskip

  Condition \ref{eqn:Admissibility} will be called the {\bf frame
  admissibility condition}.

\subsection{Real frames}\label{subsect:Real_Frames}
Let us digress momentarily in order to compare the above definitions with
their counterparts in the real case.

When the parameter $c$ is real, all the pieces of the nest intersect the real
axis. Call $I_j^n$ the intersection of $V_j^n$ with $\mathbb{R}$. The
combinatorial type of the nest is determined by how many intervals are there
left and right of $I_0^n$, the sign (orientation) of each map $g_{n,j}:I_j^n
\longrightarrow I_0^{n-1}$ and the itineraries of all $I_j^n$ through
intervals of the previous level. If we specify an arbitrary type, the {\it
unimodal admissibility conditions} are necessary so that the type can be
realized; these conditions require

\begin{itemize}
  \item Since $g_{n-1,k}$ is supposed to take $I_k^{n-1}$ onto $I_0^{n-2}$,
        the order of the intervals inside $I_k^{n-1}$ is preserved or
        reversed according to the orientation of $g_{n-1,k}$.
  \item Since $g_{n,j}:I_j^n \longrightarrow I_0^{n-1}$ is supposed to be the
        composition of all $g_{n-1,k_i}$ specified by the itinerary of
        $I_j^n$, the sign of $g_{n,j}$ must be the product of signs of the
        $g_{n-1,k_i}$ when $I_j^n$ is right of $I_0^n$ and the negative of
        that sign when $I_j^n$ is to the left of $I_0^n$ (or the other way,
        if $g_{n,0}$ reverses orientation).
\end{itemize}

\begin{figure}[h]
  \begin{center}
  \includegraphics{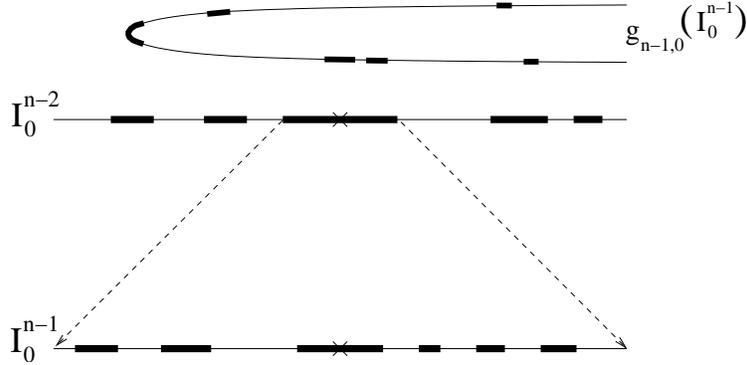}
  \caption[Unimodal admissibility conditions.]{\label{fig:Unimodality} \it
           Illustration of the unimodal admissibility conditions. The map
           $g_{n-1,0}$ spreads the intervals of level $n$ inside some
           intervals of level $n-1$. However, the order of the right
           intervals is respected and that of the left intervals is reversed.
           Note that the orientation of each left interval is also reversed
           and that $I_0^n$ maps to the leftmost position.}
  \end{center}
\end{figure}

We note first that both conditions emphasize the fact that $g_{n,0}$ is
unimodal. The first map $g_{n-1,0}$ can mix left intervals with right
intervals as in Figure \ref{fig:Unimodality}, but the order of the right
intervals is preserved and the order of the left ones is reversed (or
vice-versa). The second condition specifies that the orientation of each
$g_{n,j}$ is the product of the orientations of all intermediate steps {\it
including the fact that $g_{n-1,0}$ has different orientations on each side of
0}. The important observation to make is that the simplicity of the unimodal
admissibility conditions is due to the existence of a natural order on
$\mathbb{R}$. In the more general case of complex polynomials, the order of
intervals is replaced by relative locations of nest pieces within a frame.
The requirement that relative orders are preserved is replaced by Conditions
\ref{eqn:Admissibility} and the rule of signs is replaced by a compatible
choice of labels. \\

\subsection{Combinatorial classification}
We are ready to state the main theorem of this Section. In loose language, it
states the existence within the quadratic family, of arbitrary admissible
finite combinatorial types.

\defn{} We will say that two non-renormalizable polynomials are {\bf weakly
  combinatorially equivalent} if they have the same combinatorial types at
  every level, so that they differ only by the orientation of their frames.
\medskip

\note{} The point $g_n(0)$ is contained in $V_0^{n-1}$. In particular, it is
  possible to apply the map $g_{n-1}$ to it and, in fact, we could keep
  composing first return maps of lower levels until the first return of the
  critical orbit to $V_0^n$. This argument shows that for weak combinatorially
  equivalent maps, $g_{n+1}$ is formed by the same composition of previous
  levels first return maps and consequently, {\it the first returns to
  corresponding pieces happen at the same times}. In the next sections we will
  make use of this property.

\begin{thm}\label{thm:Frame_Renormalizations}
  Consider a finite combinatorial type of level $n$, together with a parapiece
  $\Delta$ of parameters that satisfy it up to level $n-1$. Let $\ell$ be the
  level of the last lateral return prior to level $n$ and let
  $$
    r = \left\{
          \begin{array}{ll}
            1 & \text{if } g_n \text{ is a central return} \\
            2^{n-\ell} & \text{if } g_n \text{ is lateral. }
          \end{array}
        \right.
  $$
  Then there exist $r$ parapieces inside $\Delta$ each consisting of
  parameters satisfying the same weak combinatorial type to level $n$.

  Moreover, for any such parapiece $\Delta'$, the first returns $\big\{ g_n[c]
  \mid c \in \Delta' \big\}$ form a full DH quadratic-like family.
\end{thm}

\note{} This property of accumulating powers of 2 during central cascades is
  related to the phenomenon that makes Lyubich's theorem possible. Namely, the
  fact that the moduli grow linearly from lateral return to lateral return,
  even though they decrease by half on each central return. \\

\pf{} We are already acquainted with the central symmetry of frames. It is
  obvious that the dual graph of a frame can be symmetric only about its
  critical vertex $\xi$. Because of this, the frame $F_{\ell+1}$ cannot be
  symmetric around the lateral cell $C$ where $g_{\ell}(0) \in \left(
  V_0^{\ell-1} \setminus V_0^{\ell} \right)$ falls, so the pull back
  $F_{\ell+2}$ cannot have more than 2-fold symmetry around the origin.

  By definition, the (possibly empty) sequence $\{ g_{\ell+1}, \ldots, g_{n-1}
  \}$ is the beginning of a cascade of central returns of length $n-\ell$.
  Therefore, the frame graph $\Gamma(F_{n+2})$ has exactly $(2^{n-\ell})$-fold
  symmetry around $\xi_n$. \\

  Let $c \in \Delta$. Every map $g_{n-1,k}$ takes its corresponding piece
  $V_k^{n-1}$ onto $V_0^{n-2}$. Then the pull back by $g_{n-1,k}$ of any
  region inside $V_0^{n-2}$ is well defined and located inside $V_k^{n-1}$. In
  particular, for every piece $V_j^n$ listed in the type of level $n$, the
  itinerary prescribes the sequence of returns $g_{n-1,0}, g_{n-1,k_1},
  \ldots, g_{n-1,k_r}$, so the univalent pull back of $V_0^{n-1}$ under the
  composition $(g_{n-1,k_1} \circ \ldots \circ g_{n-1,k_r})$ is a well defined
  piece inside $V_{k_1}^{n-1}$. Let us name this piece $U'_j$.

  Clearly $U'_0 \subset V_1^{n-1}$ because the itinerary of the critical piece
  $V_0^n$ begins with the first return of 0 to level $n-1$. As $c$ moves
  within $\Delta$, this return can be made to fall in $U'_0$. All $c$ with
  this property form a parapiece $\Delta^* \Subset \Delta$ that can be
  described as the set of parameters for which the itinerary of $U_0$ is as
  originally prescribed; i.e. $U_0 = V_0^n$. For the rest of the argument we
  will restrict $c$ to $\Delta^*$.

  For $j \geq 1$, the $g_{n-1,0}$-pull back of $U'_j$ will be called $U_j$;
  however, $g_{n-1,0}$ is 2 to 1, so we have to decide on a frame orientation
  before locating these pieces inside $F_{n+2}$.

  The combinatorial type of level $n$ involves the label $\sigma_{n+2}$ that
  specifies the cell in $F_{n+2}$ containing the first return $g_n(0)$. If
  this return is central there is no choice: The return falls on the piece
  $V_0^{n+1}$ inside the central cell. Otherwise, we need to recall the
  discussion above. After a (possibly vacuous) cascade of central returns,
  there are $\frac{r}2 = 2^{n-\ell-1}$ cells of $F_{n+1,k_1}$ that can be
  labeled with $\sigma_{n+2}$ and contain $U'_1$. This comes from the
  $n-\ell-1$ choices of orientation taken from level $\ell+1$ to $n-1$.
  Assuming that the return $g_n(0)$ is lateral, there is one more choice of
  orientation to make, so $F_{n+2}$ has $(2^{n-\ell})$ cells that can host
  $U_j$. Once this decision is made, the label orientation is determined and
  the rest of the pieces $U_j$ are forcibly placed around the frame $F_{n+2}$.
  \\

  We have constructed pieces $U_j \subset V_0^{n-1}$ that follow the given
  itineraries. It rests now to show that for some parameters $c \in \Delta^*$,
  the $U_j$ can be made to coincide with the respective $V_j^n$. This can be
  shown as follows. The itinerary of $V_0^n$ (and of 0) ends with the first
  return $g_n$ of 0 to $V_0^{n-1}$. This return generates a full family for $c
  \in \Delta^*$, so we can choose a parapiece $\Delta^{**}$ of $c$ such that
  $g_n(0) \in U_1$.

  The second return to $V_0^{n-1}$ is specified by the itinerary of $U_1$.
  From this observation we conclude that $U_1 = V_1^n$ from the definition of
  nest. Also, this second return generates a full family for $c \in
  \Delta^{**}$, so we can choose an even smaller parapiece $\Delta^{***}$ of
  parameters $c$ such that $g_{n,1}(0) \in U_2$. This argument can be pursued
  till the end to obtain the parapiece $\Delta'$ of values $c$ for which every
  $U_j = V_j^n$.
\QED

Repeated application of Theorem \ref{thm:Frame_Renormalizations} yields the
following.

\begin{corol}\label{corol:Infinite_Type}
  Arbitrary infinite sequences of finite, weak combinatorial types can be
  realized in the quadratic family, as long as they satisfy the admissibility
  condition at every level. The set of parameters satisfying the complete type
  is the intersection of a family of nested sequences of parapieces, with
  $2^n$ of them at every non-central level $n$.
\end{corol}

\pf{} This is clear, since each $\Delta$ contains at least one parapiece
  $\Delta'$ that satisfies the combinatorial type at level $n$. The collection
  of first return maps of level $n$ for parameters in $\Delta'$ forms a full
  family, so we can apply Theorem \ref{thm:Frame_Renormalizations} again. An
  arbitrary choice of orientation at every level gives an infinite nested
  sequence of parapieces. Evidently, a parameter in the intersection satisfies
  the prescribed combinatorics at every level.

  Every level accounts for one dyadic choice of orientation. Although they are
  not apparent during central cascades, the previous proof shows that they
  accumulate to display $2^{n-\ell}$ pieces of level $n$ inside each of the
  $2^{\ell}$ pieces of (lateral) level $\ell$.
\QED

The set of parameters that are combinatorially equivalent to a given one
cannot be completely characterized without some amount of analytical
information. Corollary \ref{corol:Infinite_Type} describes such set as a
collection of nested sequences of parametric pieces, but it does not say
whether they intersect in single points or in more complicated regions. The
fact that the parapieces shrink to a unique parameter amounts to {\it
combinatorial rigidity}; this was the strategy of Yoccoz to establish local
connectivity in the case of non-renormalizable polynomials. For such
parameters, he showed that the sum of paramoduli is infinite, so the set of
parameters in the nested intersections of parapieces becomes a Cantor set. In
particular, if the type includes no central returns, every parapiece contains
exactly two pieces of the next level and the Cantor set has a natural dyadic
structure. Thus, for some precise sequences of combinatorial types, the choice
of frame orientations at every level may single out a unique parameter.

\note{} It should be remarked that alternative classifications of
  combinatorial properties are possible and indeed quite useful. Of particular
  notice is D. Schleicher's concept of {\it internal addresses} (see
  \cite{Internal_addresses}), describing a combinatorial type in terms of an
  irreducible sequence of hyperbolic components that encodes the critical
  orbit information with increasing precision.

\section{Examples}\label{sect:Examples}
\setcounter{equation}{0}

We present here two instances of the use of our combinatorial model. Every
first renormalization type corresponds to a maximal hyperbolic component of
the Mandelbrot set; these are classified in \ref{subsect:Maximals}.  A
rotation-like map is an unimodal map whose postcritical set is semi-conjugate
to a circle rotation; the Fibonacci map being n instance. In
\ref{subsect:Rot-like} we find complex quadratic maps with the same property.
Other applications, including a classification of complex quadratic Fibonacci
polynomials, can be found in \cite{2nd_part}.

\subsection{Maximal hyperbolic components}\label{subsect:Maximals}
Consider an arbitrary combinatorial type up to some level $n$, with the
property that the last return is not central. Upon specifying a frame
orientation, there is a unique parapiece $\Delta$ consisting of parameters
that satisfy the given combinatorics. Clearly, parapieces corresponding to
different types must be disjoint.

When the return to level $n+1$ is central, there is no need to orient the
frame; that is, there is a unique piece $\Delta' \subset \Delta$ of parameters
featuring this central return. Then, if a parameter in $\Delta$ has an
infinite cascade of central returns starting at level $n+1$, its combinatorial
type will be completely determined by the initial $n$ levels. The unique
sequence of nested parapieces $\Delta \supset \Delta' \supset \ldots$
intersects in the set $M'$ of renormalizable parameters whose first $n$ nest
levels are as prescribed. It is known that $M'$ is quasi-conformally
homeomorphic to $M$ (see \cite{DH_orsay} and \cite{L_parapuzzle}). In fact,
this homeomorphism is given by {\bf straightening}: For every $c \in M'$ there
is a quasi-conformal map $h$ that realizes the conjugation $h \circ g_n = f
\circ h$ between $g_n[c]$ and some quadratic polynomial $f$; moreover, $h$
satisfies $\overline{\partial}h = 0$ on the small filled Julia set of
$g_n[c]$.

Since the parameters in $M'$ have a well defined nest, the renormalization is
not of immediate type. The base of such ``small copy'' of $M$ is a primitive
hyperbolic component $H$. Since $H$ is a quasi-conformal deformation of
$\cardioid$, its boundary has a cusp point. Also, the parameters in $H$ are
exactly once renormalizable, so $H$ is maximal (see definitions at the
beginning of Section \ref{sect:Basics}). \\

The above discussion shows that any finite frame type is associated to a
maximal hyperbolic component of $M$. Conversely, each maximal copy of $M$ is
encoded by the type of its frame, that is, by the associated graph
$\Gamma(F_{n+1})$ or its label sequence. Note that the frame graph of level
$n'>n$ consists of a bouquet of $2^{n'-n}$ copies of $\Gamma(F_{n+1})$ with
their central vertices identified. This is illustrated in the right hand
example in Figure \ref{fig:Two_Frames}. The beautiful pictures of small
Mandelbrot copies with hundreds of mini-copies spiraling in all directions
belong naturally to the class of finite nest types that conclude with a long
central cascade.

\subsection{Rotation-like maps}\label{subsect:Rot-like}
Let $c_{\text{fib}}=-1.8705286321 \ldots$ parametrize the Fibonacci map $z
\mapsto z^2+c_{\text{fib}}$. This is the unique real quadratic polynomial with
the property that the critical orbit has closest returns to 0 exactly when the
iterates are the Fibonacci numbers; see \cite{LM_fibo}. In terms of the
principal nest, $f_{c_{\text{fib}}}$ satisfies the equivalent
condition:\vspace{.15in}

\begin{center}\begin{minipage}{5.5in}
For $n \geq 2$, each level of the principal nest consists of the central piece
$V_0^n$ and a unique lateral piece $V_1^n$. The first return map of previous
level $g_{n-1}:V_0^{n-1} \longrightarrow V_0^{n-2}$ interchanges the central
and lateral roles:
$$g_{n-1}(V_0^n) \Subset V_1^{n-1}, g_{n-1}(V_1^n) = V_0^{n-1}.$$

Additionally, the first returns to $Y_0^{(1)}$ and $V_0^0$ happen on the third
and fifth iterates respectively.
\end{minipage}\end{center}\vspace{.15in}

To discern the critical orbit behavior of $f_{c_{\text{fib}}}$, note that
every level of the nest has a unique lateral piece and so, in a sense, every
first return comes as close as possible to being central without actually
being so.  This means that the map $f_{c_{\text{fib}}}$ is not renormalizable
in the classical sense, although its combinatorics can be described as an
infinite cascade of {\it Fibonacci renormalizations} in the space of {\bf gql}
maps with one lateral piece.

The Fibonacci map features as a decisive case in the proof of Lyubich's
theorem; see \cite{L_nest}. Here we will describe a family of unimodal maps
with similar behavior and extend it to a family of complex quadratic maps. \\

Let $S = (S_0, S_1, \ldots)$ be a strictly increasing sequence of numbers
such that $\frac{S_{j+1}}{S_j} \leq 2.$ The {\it $S$-odometer} is a symbolic
dynamical system $(\Omega,T)$ defined as follows. For any nonnegative $n$
there is a $k$ such that $S_k \leq n < S_{k+1}.$ Then $n = S_k + n_1$ with
$n_1 < S_k.$ By splitting further $n_1 = S_{k'} + n_2$ (with $k'<k$ and
$n_2 < S_{k'}$) and so on, we obtain the decomposition
$$n = d_k \cdot S_k + \ldots + d_0 \cdot S_0$$
where each $d_j$ is either 0 or 1. Letting $d_j =0$ for $j>k,$ we get the
sequence
$$\langle n \rangle = (d_0, d_1, \ldots) \in \{ 0,1 \}^{\mathbb{N}}.$$

We use $\langle \mathbb{N} \rangle$ to denote
$\{ \langle n \rangle \mid n \in \mathbb{N} \}$ and let $\Omega$ be the
closure
$$
  \Omega = \overline{\langle \mathbb{N} \rangle} =
  \{ \omega \in \{ 0,1 \}^{\mathbb{N}} \mid
  \sum_{i=0}^j \omega_j S_i < S_{j+1} \text{ for all } j \geq 0 \}.
$$

The map
$T:\langle \mathbb{N} \rangle \longrightarrow \langle \mathbb{N} \rangle$
is given by $T \langle n \rangle = \langle n+1 \rangle.$ This map does not
always extend uniquely to $\Omega.$ When there is an extension, the
dynamical system $(\Omega,T)$ obtained from the sequence $S$ is called a
{\bf $S$-odometer}. It can be described as an adding machine with variable
stepsize. \\

Let us relate the above concept to interval dynamics. First, some definitions.

Consider a unimodal map $f:I \longrightarrow I$ where $I = [c_1,c_2]$ and
$\{ 0, c_1, c_2, \ldots \}$ is the critical orbit. Let $D_1 =[c_1,0]$ and, for
$n \geq 2,$ define
$$
  D_{n+1} =
  \left\{
  \begin{array}{ll}
    [ c_{n+1},c_1 ] & 0 \in D_n \\
    f(D_n)          & 0 \nin D_n
  \end{array}
  \right.
$$

The sequence $S = (S_0, S_1, \ldots )$ of {\it cutting times} consists of
those $n$ such that $0 \in D_n.$ Note that $S_0 = 1.$ It is easy to show that
$S_{k+1}-S_k$ is also a cutting time so we can define the {\bf kneading map}
$Q:\mathbb{N} \longrightarrow \mathbb{N}$ by the relation
$$S_{Q(k)} = S_{k+1} - S_k.$$

\begin{lemma}
  If $S$ is the sequence of cutting times of a unimodal map $f,$ the following
  characterization of $\Omega$ holds:
  $$
    \Omega = \{ \omega \in \{ 0,1 \}^{\mathbb{N}} \mid
    \omega_j=1 \Rightarrow \omega_i = 0 \text{ for } Q(j+1) \leq i \leq j-1 \}.
  $$
  Also, if $Q(k) \longrightarrow \infty,$ then $T$ extends uniquely to
  $\Omega$ and is conjugate to the action of $f$ on its postcritical set.
\end{lemma}

See \cite{B_odometers} for proofs. \\

In the case of the Fibonacci polynomial, the above definitions correspond to
the description of the critical orbit in Subsection 3 of \cite{LM_fibo}. There
it is shown that $(\Omega,T)_{c_{\text{fib}}}$ is semiconjugate to the
circle rotation by $\rho = \frac{\sqrt{5}-1}2.$ Real {\bf rotation-like
maps}, as defined in \cite{B_odometers}, are unimodal maps that generalize
this behavior.

Let $\rho \in [0,1) \setminus \mathbb{Q}$ with continued fraction expansion
$\rho = [a_1,a_2,\ldots]$ and denote its convergents with $\frac{p_i}{q_i}$
so that $\frac{p_0}{q_0} = \frac01$ and $\frac{p_1}{q_1} = \frac1{a_1}.$

\begin{thm}{\rm \cite{B_odometers}}
  Consider the sequence $r_k$ starting with $r_1 = q_1-1$ and whose
  $(k+1)^{st}$ element is given recursively by $r_{k+1} = r_k + a_{k+1}.$ Then
  the $S$-sequence given by
  $$
    \begin{array}{rcll}
      S_{r_k} & = & q_k \\
      S_{r_k+j} & = & (j+1)q_k & \text{ for } 1 \leq j < a_{k+1}
    \end{array}
  $$
  is realized as the sequence of cutting times of some quadratic polynomial.
  Moreover, the application
  $$
    \Pi_{\rho}(\omega) = \sum \omega_j S_j \rho \, (\text{mod} \, 1)
  $$
  from $\Omega$ to the unit circle is well defined and continuous. This map
  satisfies $\Pi_{\rho} \circ T = R_{\rho} \circ \Pi_{\rho},$ where
  $R_{\rho}$ is the rotation by angle $\rho,$ and is 1 to 1 everywhere except
  at the preimages of 0.
\end{thm}

In terms of the principal nest, the behavior that characterizes rotation-like
maps is a succession of central cascades followed by one lateral escape. That
is, the critical orbit falls in $V_0^{S_k-1}$ starting a central cascade.
After iterating the first return map $g_k$ for $a_k-1$ turns, we get a lateral
return on $V_1^{S_k}.$ Next, $g_{S_k,1}$ creates a new cascade and so on. In
particular, the Fibonacci map is the special case of a rotation-like map where
every central cascade has length 0. \\

Consider an arbitrary sequence $a_1, a_2, \ldots$ of positive integers. We
will construct now a Cantor set of {\it complex} rotation-like parameters with
central cascades of length $a_i-1.$ By theorem
\ref{thm:Frame_Renormalizations}, it is only necessary to give an admissible
description of labeling sequences and to show that it models the combinatorics
mentioned above.

The initial labeling data for our map is $q=2$ and $\sigma_0 =$ {\sf '1'}, so
rotation-like maps will all be located in the $1/2$-limb. Note also that on
central return levels, $\sigma_{k+1} = {\sf '0'}\sigma_k.$ Therefore, we only
need to specify the labels $\sigma_{r_\eta}$ for $r_\eta = \sum a_j.$

Let $(\tau_1,\tau_2,\ldots)$ be a sequence of random chains of {\sf 'l'}'s and
{\sf 'r'}'s so that $\tau_i$ has length $a_i+1.$ Set $\sigma_{r_1} = \tau_1
{\sf '0'}$ and $\sigma_{r_2} = \tau_2 \sigma_{r_1-1} = \tau_2 {\sf '00 \ldots
01'}.$ Now we can define inductively $\sigma_{r_j} = \tau_j {\sf '0'}
\sigma_{r_{j-1}-1}.$

\begin{prop}
  The label sequence $(q;\sigma_0, \sigma_1, \ldots )$ defined above is
  admissible, it completely describes a combinatorial type and the
  corresponding map is rotation-like.
\end{prop}

\pf{}
  The fact that the sequence of labels determines the type can be seen to be
  true since there are no consecutive lateral returns. This implies that the
  nest has exactly one lateral piece at those levels (and none elsewhere) so
  its position within the frame is completely determined by $\sigma_{r_j}.$

  As mentioned above, ${\sf '0'}\sigma_k$ (when $k \neq r_j$) is an admissible
  label since it corresponds to the central cell of $F_{k+1}.$ Now consider
  what happens to the central cell labeled ${\sf '0'} \sigma_{r_{j-1}-1}.$
  Since level $r_{j-1}$ corresponds to a non-central return, $F_{r_{j-1}+1}$
  has two preimages of that cell, labeled ${\sf 'l0'} \sigma_{r_{j-1}-1}$ and
  ${\sf 'r0'} \sigma_{r_{j-1}-1}$ respectively. On consecutive central
  returns, we double the number of pull-backs of such cells and thus, use all
  possible combinations of ${\sf 'l'}$ and ${\sf 'r'}$ to label them. A glance
  to the frame graph shows that these are the cells neighboring the central
  one (see \cite{Daniel}). An eventual lateral return must fall precisely in
  one of these cells, and this is what happens when $\sigma_{r_j} = \tau_j
  {\sf '0'} \sigma_{r_{j-1}-1}.$
\QED

The real rotation-like maps studied in \cite{B_odometers} correspond to a
careful choice of the $\tau_j$. In fact, it is possible to extract a kneading
sequence from the rotation number data. Then, a result of Yoccoz guarantees
that there is a unique real polynomial in that combinatorial class.

The complex maps corresponding to other choices of $\tau_j$'s have the same
weak combinatorial behavior, so the critical orbits of two maps with the same
sequence $a_1,a_2,\ldots$ are conjugate. In particular we obtain the following
result.

\begin{corol}
  Given the sequence $a_1, a_2, \ldots$ there exists an infinite family of
  complex quadratic polynomials for which the postcritical set is conjugate to
  an $S$-odometer and semi-conjugate to the circle rotation of angle $\rho =
  [a_1, a_2, \ldots].$
\end{corol}

\section*{Appendix}
\setcounter{subsection}{0} \setcounter{thm}{0}
\renewcommand{\thesubsection}{\Alph{subsection}}
\renewcommand{\thethm}{\thesubsection.\arabic{thm}}

\subsection{Holomorphic motions of puzzle pieces and winding number}
Consider the following

\defn{} Let $X_* \subset \overline{\mathbb{C}}$ be an arbitrary set and
  $\Delta \subset \mathbb{C}$ a simply connected domain with $*$ as a base
  point. A {\bf holomorphic motion} of $X_*$ over $\Delta$ is a family of
  injections $h_{\lambda}:X_* \longrightarrow \overline{\mathbb{C}}$ $(\lambda
  \in \Delta)$ such that for each fixed $x \in X_*$, $h_{\lambda}(x)$ is a
  holomorphic function of $\lambda$ and $h_* = \text{id}$.  For every $\lambda
  \in \Delta$ we write $X_{\lambda}$ to denote the set $h_{\lambda}(X_*)$.
\medskip

Holomorphic motions are extremely versatile owing to their regularity
properties. The motion can always be extended beyond $X_*$ and is
transversally quasi-conformal. This is the content of the {\it
$\lambda$-lemma}.

\begin{thm}{\rm \cite{Final_lambda}, \cite{Lambda_lemma}
            ({\bf the $\lambda$-lemma})}
  For every holomorphic motion
  $h_{\lambda}:X_* \longrightarrow \overline{\mathbb{C}}$, there is an
  extension to a holomorphic motion
  $H_{\lambda}:\overline{\mathbb{C}} \longrightarrow \overline{\mathbb{C}}$.
  The extension to the closure $\overline{h}_{\lambda}:\overline{X}_*
   \longrightarrow \overline{\mathbb{C}}$ is unique. Moreover, there is a
  function $K(r)$ approaching 1 as $r \rightarrow 0$ such that the maps
  $h_{\lambda}$ are $K(r)$-quasi-conformal, where
  $r = \text{d}_{\Delta}(*,\lambda)$ is the hyperbolic distance between $*$
  and $\lambda$ in $\Delta$.
\end{thm}

We are interested in the case when the holomorphic motion is defined over a
parapiece $\Delta$ of $M$. In agreement with the notation used in the main
body of this work, we use $c$ instead of the classical $\lambda$ to denote
parameters in $\Delta$. When an object is defined for any $c \in \Delta$, we
express its dependence on the parameter by writing $\text{OBJ}[c]$.

As mentioned in Section \ref{sect:Basics}, $\Delta$ can be interpreted as the
set of parameters for which {\it a given combinatorial behavior holds, up to a
return $g(0)$ of the critical orbit to some central piece $V$}. In particular,
this description provides a natural base point for $\Delta$. Namely, the
superattracting parameter $c_0$ for which $g_{c_0}(0) = 0$. The little
$M$-copy associated to $\Delta$ can be defined as the set of parameters for
which the iterates $\{ g_c(0), g_c^{\circ 2} (0), \ldots \}$ remain in $V[c]$
(refer to Subsection \ref{subsect:Maximals}).

The dynamics in the region $N_c$ (defined at the beginning of Subsection
\ref{subsect:External_Rays_n_Wakes}) is always conjugate to $z \mapsto z^2$,
so varying the parameter $c \in \mathbb{C}$ provides a holomorphic motion of
any specified (open) ray or equipotential. When $c$ is restricted to $\Delta$,
the combinatorics require that some rays land together, enclosing the boundary
of $V[c]$. Since the intersection $\bdry V \cap K$ is a collection of
preimages of the fixed point $\alpha$ and these vary holomorphically with $c$,
there is a natural holomorphic motion of $\bdry V[c_0]$ over $\Delta$. This
can be extended to a holomorphic motion $h_c:V[c_0] \longrightarrow V[c]$.

The holomorphic motion of a puzzle piece can be viewed as a complex
1-dimensional foliation of the bi-disk
$$\mathbb{V} = \bigcup_{c \in \Delta} V[c] \in \mathbb{C}^2$$
whose leaves are the graphs of the functions $c \mapsto h_c(p)$ for every
$p \in V[c_0]$. Under this interpretation we will write
$\{ c \mapsto V[c] \mid c \in \Delta \}$ to refer to the motion.

\defn{}
  A correspondence $c \mapsto \phi(c)$ such that $\phi(c) \in V[c]$ determines
  a section $\phi:\Delta \longrightarrow \mathbb{V}$ of the holomorphic motion
  $h$. It is said to be a {\bf proper holomorphic section} if it maps
  $\bdry \Delta$ into the torus
  $\delta\mathbb{V} = \bigcup_{c \in \bdry \Delta} \bdry V[c]$.

  We say that a proper section $\{ c \mapsto \phi(c) \}$ has {\bf winding
  number} $n$ if the curve $\phi(\bdry \Delta)$ has winding number $n$ with
  respect to the vertical generator of the 1-dimensional homology of
  $\delta\mathbb{V}$.
\medskip

In the case $\phi(c) = g_c(0)$, this return map determines a proper section
since $ g_c(0) \in V[c]$ for all $c$ and $c \in \bdry \Delta \Rightarrow
g_c(0) \in \bdry V[c]$. Each return map $g_c:g_c^{-1}(V) \longrightarrow V$ is
a quadratic-like map and the associated map
$$g_c:\mathbb{U} \longrightarrow \mathbb{V},$$
where $\mathbb{U} = \bigcup g_c^{-1}(V[c])$, is called a {\bf DH
quadratic-like family}. We can interpret intuitively the fact that a family
has winding number $n$ as saying that, as $c$ goes once along $\bdry \Delta$,
the point $g_c(0)$ goes $n$ times around the (moving) boundary of the piece
$V[c]$.

An immediate consequence of extending the holomorphic motion of $\bdry
V[c_0]$, is the fact that $\{ g_c \mid c \in \Delta \}$ is a full family; that
is, there is a homeomorphism $\text{Hyb}:\Tilde{M} \longrightarrow \Delta$
from a neighborhood $\Tilde{M}$ of $M$ to $\Delta$ with the following
property: For every parameter $c' \in \Tilde{M}$, $g_{\text{Hyb}(c')}$ is
hybrid equivalent\footnote{see Subsection \ref{subsect:Maximals}} to $z
\mapsto z^2+c'$.  This of course, justifies the existence of the small
$M$-copy associated to $\Delta$.


\begin{thebibliography}{99999}


\bibitem[BH1]{BH_cubics1} B. Branner, J. Hubbard.
  {\it The iteration of cubic polynomials. Part I: The global topology of
  parameter space.}
  Acta Math., {\bf 160} (1988), 143-206.


\bibitem[BH2]{BH_cubics2} B. Branner, J. Hubbard.
  {\it The iteration of cubic polynomials. Part II: Patterns and
  parapatterns.}
  Acta Math., {\bf 169} (1992), 229-325.


\bibitem[BKP]{B_odometers} H. Bruin, G. Keller and M. st. Pierre.
  {\it Adding machines and wild attractors.}
  Ergod. Th. \& Dynam. Sys., {\bf 17} (1997), 1267-1287.


\bibitem[D2]{Chirurgie} A. Douady.
  {\it Chirurgie sur les applications holomorphes.}
  In: Proc. ICM, Berkeley, (1986), 724-738.


\bibitem[DH1]{DH_orsay} A. Douady \& J. H. Hubbard,
  {\it \'Etude dynamique des polyn\^omes complexes I \& II.}
  Publ. Math. Orsay, 1984-85.


\bibitem[DH2]{DH_p-l} A. Douady \& J. H. Hubbard,
  {\it On the dynamics of polynomial-like maps.}
  Ann. Sci. \'Ec. Norm. Sup., {\bf 18} (1985), 287-343.


\bibitem[GLT]{GLT_from_B} P. J. Grabner, P. Liardet and R. F. Tichy.
 {\it Odometers and systems of enumeration.}
 Acta Arithmetica, {\bf 70} (1995), 103-123.


\bibitem[H]{Tableaux} J. H. Hubbard,
  {\it Local connectivity of Julia sets and bifurcation loci: Three theorems of
  J.-C. Yoccoz.}
  In: Topological Methods in Modern Mathematics pp. 467-511 (ed. L. Goldberg
  \& A. Phillips), (Publish or Perish, 1993).


\bibitem[J]{J} W.Jung.
  PC software {\tt mandel.exe}; available at: {\tt
  http://www.iram.rwth-aachen.de/$\sim$jung/indexp.html}


\bibitem[LS]{Internal_addresses} E. Lau and D. Schleicher.
  {\it  Internal Addresses of the Mandelbrot Set and Irreducibility of
  Polynomials.}
  Preprint IMS at Stony Brook, \# 1994/19. 


\bibitem[L1]{L_attractor} M. Lyubich.
  {\it Combinatorics, geometry and attractors of quasi-quadratic maps.}
  Ann. of Math., {\bf 140}, (1994), 347-404.


\bibitem[L2]{L_nest} M. Lyubich.
  {\it Dynamics of quadratic polynomials, I-II.}
  Acta Math., {\bf 178} (1997), 185-297. 


\bibitem[L3]{L_parapuzzle} M. Lyubich.
  {\it Dynamics of quadratic polynomials, III. Parapuzzle and SBR measures.}
  In: G\'eom\'etrie Complexe et Syst\'emes Dynamiques. Volume in Honor of
  Adrien Douady's 60th Birthday. Ast\'erisque 261, (2000), 173-200.


\bibitem[LM]{LM_fibo} M. Lyubich and J. Milnor.
  {\it The Fibonacci unimodal map.}
  J. Amer. Math Soc., {\bf 6} (1993), 425-457.


\bibitem[MSS]{Lambda_lemma} R. Ma\~n\'e, P. Sad and D. Sullivan.
 {\it On the dynamics of rational maps.}
 Ann. Sci. \'Ec. Norm. Sup., {\bf 4}, (1983), 193-217.


\bibitem[Ma]{Martens} M. Martens.
  {\it Distortion results and invariant Cantor sets of unimodal maps.}
  Ergod. Th. \& Dynam. Sys., {\bf 14}, (1994), 331-349.


\bibitem[M1]{M_book} J. W. Milnor,
  Dynamics in One Complex Variable.
  (Vieweg, 1999).


\bibitem[M2]{Portraits} J. W. Milnor,
  {\it Periodic orbits, external rays and the Mandelbrot set: An expository
  account.}
  In: Asterisque 261 `Geometrie Complexe et Systemes Dynamiques', pp.
  277-333, (SMF 2000).


\bibitem[P]{2nd_part} R. P\'erez,
  {\it Geometry of $Q$-recurrent maps.}
  In preparation.

\bibitem[R]{Roesch} P. Roesch.
  {\it Holomorphic motions and puzzles (following Shishikura).}
  In: The Mandelbrot set, theme and variations, edited by Tan Lei, LNS 274,
  Cambridge, (2000).


\bibitem[Sl]{Final_lambda} Z. Slodkowsky.
  {\it Holomorphic motions and polynomial hulls.}
  Proc. Amer. Math. Soc., {\bf 111}, (1991), 347-355.


\bibitem[Sm]{Daniel} D. Smania.
  {\it Puzzle geometry and rigidity: The Fibonacci cycle is hyperbolic.}
  {\bf arXiv}: {\tt math.DS/0203164.}


\end{thebibliography}
\end{document}